\documentclass{amsart}
\usepackage{enumerate, amssymb}

\def\pdfsyncstart{}
\def\pdfsyncstop{}

\def\bdi{\pdfsyncstop\begin{diagram}}
\def\edi{\end{diagram}\pdfsyncstart}

\usepackage{diagrams}
\diagramstyle[scriptlabels,small]

\newarrow{Equal}{}{=}{}{=}{}

\def\lto{{\longrightarrow}}
\def\hto{{\hookrightarrow}}
\def\xto{\xrightarrow}

\newcommand{\ff}{{\mathfrak f}}
\newcommand{\fg}{{\mathfrak g}}
\newcommand{\fh}{{\mathfrak h}}

\newcommand{\fX}{{\mathfrak X}}
\newcommand{\fY}{{\mathfrak Y}}
\newcommand{\fZ}{{\mathfrak Z}}

\newcommand{\cala}{{\mathcal A}}
\newcommand{\calb}{{\mathcal B}}
\newcommand{\calc}{{\mathcal C}}
\newcommand{\cald}{{\mathcal D}}

\newcommand{\calf}{{\mathcal F}}
\newcommand{\calg}{{\mathcal G}}
\newcommand{\calh}{{\mathcal H}}

\newcommand{\calm}{{\mathcal M}}
\newcommand{\caln}{{\mathcal N}}
\newcommand{\calo}{{\mathcal O}}
\newcommand{\calp}{{\mathcal P}}
\newcommand{\calq}{{\mathcal Q}}
\newcommand{\calr}{{\mathcal R}}
\newcommand{\cals}{{\mathcal S}}
\newcommand{\calt}{{\mathcal T}}

\newcommand{\calow}{{\calo_W}}
\newcommand{\calox}{{\calo_X}}
\newcommand{\caloy}{{\calo_Y}}

\def\calas{{\cala_*}}
\def\calbs{{\calb_*}}    
\def\calds{{\cald_*}}

\def\calms{{\calm_*}}
\def\calns{{\caln_*}}
\def\calps{{\calp_*}}
\def\calqs{{\calq_*}}
\def\calrs{{\calr_*}}
\def\calss{{\cals_*}}
\def\calts{{\calt_*}}
\def\caloxs{{\calo_{X_*}}}
\def\calows{{\calo_{W_*}}}

\def\cExt{{\mathcal E}xt}
\def\cTor{{\mathcal T}or}

\def\cHH{{\mathcal H\mathcal H}}

\newcommand{\bbbb}{{\mathbb B}}
\newcommand{\bbbh}{{\mathbb H}}
\newcommand{\bbbl}{{\mathbb L}}
\newcommand{\bbbq}{{\mathbb Q}}
\newcommand{\bbbc}{{\mathbb C}}

\newcommand{\bbbs}{{\mathbb S}}

\newcommand{\bbbp}{{\mathbb P}}

\newcommand{\bbbz}{{\mathbb Z}}

\newcommand{\bdot}{\sbullet}

\newcommand{\vp}{\varphi}

\DeclareMathOperator{\dotimes}{\underline{\otimes}}

\DeclareMathOperator{\Ext}{Ext}

\DeclareMathOperator{\HH}{HH}

\DeclareMathOperator{\Hom}{Hom}
\DeclareMathOperator{\id}{id}

\DeclareMathOperator{\Spec}{Spec}

\DeclareMathOperator{\Tor}{Tor}

\DeclareMathOperator{\Mor}{Mor}

\newcommand{\sbullet}{{\scriptscriptstyle\bullet}}

\DeclareMathOperator{\Mod}{\ensuremath{ \mathbf{Mod}}}

\newcommand{\cata}{\ensuremath{ \mathbf A}}
\newcommand{\catc}{\ensuremath{ \mathbf C}}

\newcommand{\catf}{\ensuremath{ \mathbf F}}
\newcommand{\catm}{\ensuremath{ \mathbf M}}
\def\catmor{{\mathbf{Mor}}}

\def\Cech{\v{C}ech\ }

\def\calcb{{\calc^\sbullet}}

\theoremstyle{definition}
\newtheorem{defn}{Definition}[subsection]

\theoremstyle{plain}
\newtheorem{prop}[defn]{Proposition}
\newtheorem{theorem}[defn]{Theorem}
\newtheorem{lem}[defn]{Lemma}
\newtheorem{lemma}[defn]{Lemma}
\newtheorem{cor}[defn]{Corollary}

\theoremstyle{remark}
\newtheorem{rem}[defn]{Remark}
\newtheorem{rems}[defn]{Remarks}

\newtheorem{sit}[defn]{}
\newtheorem{exam}[defn]{Example}

\def\bpro{\begin{prop}}
\def\epro{\end{prop}}

\def\bthm{\begin{theorem}}
\def\ethm{\end{theorem}}

\def\bdfn{\begin{defn}}
\def\edfn{\end{defn}}

\def\brem{\begin{rem}}
\def\erem{\end{rem}}

\def\bsit{\begin{sit}}
\def\esit{\end{sit}}

\def\blem{\begin{lem}}
\def\elem{\end{lem}}

\def\ba{\begin{array}}
\def\ea{\end{array}}

\def\bnum{\begin{enumerate}}
\def\enum{\end{enumerate}}

\def\be{\begin{equation}}
\def\ee{\end{equation}}

\def\bproof{\begin{proof}}
\def\eproof{\end{proof}}

\begin{document}
\title[Global Hochschild (co-)homology]{Global Hochschild (co-)homology of singular spaces}

\author{Ragnar-Olaf Buchweitz}
\address{Dept.\ of Computer and Mathematical Sciences, University of
Toronto at Scarborough, Toronto, Ont.\ M1C 1A4, Canada}
\email{ragnar@math.toronto.edu}

\author{Hubert Flenner}
\address{Fakult\"at f\"ur Mathematik der Ruhr-Universit\"at,
Universit\"atsstr.\ 150, Geb.\ NA 2/72, 44780 Bochum, Germany}
\email{Hubert.Flenner@rub.de}

\thanks{The authors were partly supported by NSERC grant
3-642-114-80, by the DFG Schwerpunkt ``Global Methods in Complex
Geometry", and by the Research in Pairs Program at
Math.~Forschungsinstitut Oberwolfach funded by the
Volkswagen Foundation.}


\begin{abstract}
We introduce Hochschild (co-)homology of morphisms of schemes or
analytic spaces and study its fundamental properties. In analogy with the
cotangent complex we introduce the so called (derived) Hochschild 
complex of a morphism; the Hochschild cohomology and homology 
groups are then the $\Ext$ and $\Tor$ groups of that complex. 
We prove that these objects are well defined, extend the known cases, 
and have the expected functorial and homological properties such 
as graded commutativity of Hochschild cohomology and existence 
of the characteristic homomorphism from Hochschild cohomology to 
the (graded) centre of the derived category.
\end{abstract}

\maketitle

{\footnotesize\tableofcontents}

\section*{Introduction}

The aim of this note is to define Hochschild (co-)homology in the global setting, 
for morphisms of schemes or analytic spaces.
 
Hochschild homology for a {\em flat\/} morphism of any type of spaces $X\to Y$ should 
coincide with $\Tor_{\bdot}^{X\times_Y X}(\calox,\calox)$, and if $Y$ is a simple point, Hochschild cohomology should agree with $\Ext^{\bdot}_{X\times X}(\calox,\calox)$, 
where $X$ is considered as a subspace of $X\times_Y X$ via the diagonal embedding. 
In the algebraic, case, when $X$ is a quasi-projective scheme over some field 
$K$, Swan \cite{Sw} showed that this requirement holds for any
globalization of the concept of Hochschild {(co-)\-homology} that had been 
proposed earlier, e.g. in \cite{GS2}. He also proved that for $X$ smooth the corresponding 
Hodge spectral sequences agree, thus the Hochschild-Kostant-Rosenberg (HKR) decomposition that results from the degeneration of that spectral sequence 
for $X$ smooth over a field $K$ of characteristic zero is the same in those 
globalizations. 

To go beyond the flat or absolute situation, with the aim to include both 
non-flat morphisms of schemes and the complex analytic case, let us first 
review the pertinent issues in the affine case. For homomorphisms $A\to B$ of 
commutative rings, even more generally, for any associative $A$--algebra $B$, 
the bar complex $\bbbb$ provides a canonical resolution of $B$ as a (right) module over the enveloping algebra $B^e:=B^{op}\otimes_AB$, and classically Hochschild (co-)homology is the (co-)homology of (the dual of) that complex, that is,
\begin{align*}
\HH_{\bdot}^{B/A}(M) := H_{\bdot}(\bbbb\otimes_{B^{e}}M)\quad\text{and}\quad \HH^{\bdot}_{B/A}(M) := H^{\bdot}(\Hom_{B^{e}}(B,M))\,,
\end{align*}
for any $B$--bimodule $M$. If $\bbbp$ is a projective resolution of $B$ as a $B^{e}$--module, then there exists a comparison map of complexes $\bbbp\to \bbbb$ over the identity on $B$, and that map is unique up to homotopy. Consequently, there are natural comparison maps
\begin{align*}
\alpha_{M}:\Tor_{\bdot}^{B^{e}}(B,M)\lto \HH_{\bdot}^{B/A}(M) 
\quad\text{and}\quad
\alpha^{M}: \HH^{\bdot}_{B/A}(M)\to \Ext^{\bdot}_{B^{e}}(B,M)\,.
\end{align*}
The maps $\alpha_{M}$ are isomorphisms whenever $B$ is flat over $A$, while the
maps $\alpha^{M}$ are isomorphisms as soon as $B$ is projective as an $A$--module.

For a morphism of analytic algebras,  one can, in principle,  mimic this approach, 
replacing the standard tensor product by the analytic one, to obtain an analytic bar complex and analytic Hochschild (co-)homology.
However, already in simple situations, for instance for a free power
series ring extension, this construction will not return the expected cohomology
groups. The main reason for this is that even for a free analytic
extension $A\to B$ of analytic algebras, $B$ is not a projective module 
over $A$ unless $A$ is artinian; see, for example, \cite{BFl2} and \cite[Satz 9]{Wol}.

The same caveat applies to the obvious extension of this definition to spaces:
in any category of ringed spaces that admits fibre products, one can form, for a morphism
$f:X\to Y$ and each integer $n\ge 1$,  the $n$-fold fibre product $X^{\times_{Y}n}:= X{\times_{Y}\cdots\times_{Y}X}$, and then restrict the structure sheaf topologically to the diagonal $\Delta_{n}:X\subseteq X^{\times_{Y}n}$. With the usual definition of the differential, one obtains the {\em sheafified (analytic) bar complex\/} on $X$, a complex of $\Delta_{2}^{-1}\calo_{X\times_{Y}X}$--modules, that agrees with the complex of the same name considered by Swan in the algebraic setting.

Instead we follow the approach introduced by Quillen in \cite{Qui3} for
morphisms of algebras. To deal with the case when $B$ is not $A$-flat, 
one replaces $B^e$ by the derived tensor product $B^{op}\dotimes_AB$ 
that can be defined as a real world algebra using free simplicial resolutions 
of $B$, and then the same definition as before applies. 
There is a natural homomorphism 
$B^{op}\dotimes_{A}B\to B^{op}\otimes_{A}B=B^{e}$ over the identity of $B$, and this 
implies the existence of natural comparison maps
\begin{align*}
\beta_{M}&:\Tor_{\bdot}^{B^{op}\dotimes_{A}B}(B,M)\lto \Tor_{\bdot}^{B^{e}}(B,M) 
\quad\text{and}\quad\\
\beta^{M}&: \Ext^{\bdot}_{B^{e}}(B,M)\lto  \Ext^{\bdot}_{B^{op}\dotimes_{A}B}(B,M)\,.
\end{align*}
The maps $\beta_{M}$ are again isomorphisms for $B$ flat over $A$, while
the maps $\beta^{M}$ are isomorphisms as long as $\Tor_{i}^{A}(B,B)=0$ for $i\neq 0$,
in particular for $B$ flat over $A$.

One advantage of this setup, as shown in  \cite{Qui3},
is that the crucial Hochschild-Kostant-Rosenberg (HKR) decomposition 
theorem for Hochschild (co-)homology generalizes to arbitrary morphisms between
commutative rings of characteristic zero, with the module of K\"ahler differential forms replaced by the cotangent complex of $A\to B$. Moreover, the elegant treatment of 
the Eckmann-Hilton argument by Suarez-Alvarez in \cite{SA} yields essentially 
automatically that the natural ring structure on $\Ext^{\bdot}_{B^{op}\dotimes_{A}B}(B,B))$ is graded commutative, thus providing the counterpart to Gerstenhaber's fundamental result in \cite{Ger} for classical Hochschild cohomology.

The second advantage is that, following Quillen's guidance, one may extend the 
technique to not necessarily flat morphisms $X\to Y$ between schemes or 
analytic spaces by replacing $X\times_Y X$ with the derived fibre product. 
We describe the more complicated case of complex analytic algebras, 
asking the reader to make the simplifications that occur in the case of 
schemes. Locally, one proceeds as follows.
Given a morphism of analytic algebras $A\to B$, resolve 
$B$ by a free DG $A$-algebra $R$, so that $R$ is a DG algebra 
with an $A$-linear derivation $\partial$ as differential\footnote{
We use throughout the convention that differentials {\em increase\/} 
degrees by 1.} and comes 
equipped with a quasiisomorphism $R\to B$ over $A$. 
Here ``free" means that $R$ is concentrated in degrees $\le 0$, that
$R^0$ is a free analytic power series ring over $A$, and that 
$R$ is free as a graded $R^0$-algebra in the usual sense. 
The derived analytic tensor product $B\tilde\dotimes_AB$
is then represented by the analytic tensor product 
$S:=R\tilde\otimes_AR 
:= R\otimes_{R^{0}}(R^{0}\tilde\otimes_AR^{0})\otimes_{R^{0}}R$
that inherits naturally the structure of a free DG algebra over $A$. The
multiplication map $S\to R$ followed by the quasiisomorphism $R\to B$
endows $B$ with a DG-module structure over $S$.
Accordingly we will define the Hochschild complex
$\bbbh_{B/A}$ as the derived tensor product $B\dotimes_SB$, 
a complex  in the derived category,
with the Hochschild (co-)\-homology functors given by  
$\Tor_{\bdot}^B(\bbbh_{B/A},- )$ and $\Ext^{\bdot}_B(\bbbh_{B/A},-)$. Clearly,
this definition works locally on any analytic space.

We should point out here that contrary to the theory of the cotangent complex that is
quite sensitive to the characteristic and requires simplicial instead of DG algebra resolutions
in case of positive or mixed characteristic, the construction of the Hochschild complex is oblivious to the characteristics of the rings involved, whence one can use DG algebras throughout.

To globalize, we use {\em resolvents\/}  as developed in 
\cite{Pal, Fle1, BFl}. This allows to define first a Hochschild complex 
on the simplicial space over $Y$ that is associated to a locally finite 
covering of the given space $X$ by Stein compact sets, using the local
construction and propagating it along the nerf of the simplicial structure, 
and then to descend to $X$ via a \Cech construction to obtain a Hochschild complex
$\bbbh_{X/Y}$ in the derived category $D^{-}(X)$. 
We verify that this procedure indeed leads to a notion of 
Hochschild (co-)homology with the expected properties, such as graded
commutativity of the natural product on Hochschild cohomology and existence of a characteristic ring homomorphism from the Hochschild cohomology to the graded centre of the derived category. For a flat morphism, this Hochschild complex reduces to
the derived tensor product 
$\calox\dotimes_{\calox\tilde\otimes_{\caloy}\calox}\calox 
=(L\Delta^{*})\Delta_{*}\calox\in D^{-}(X)$, where $\Delta: X\to X\times_{Y}X$ is the 
diagonal embedding, and so the requirements laid out at the beginning are satisfied.
 
In a subsequent paper we will show that the 
HKR-decomposition theorem holds as well, thus globalizing Quillen's theorem. 
We note that in the meantime F.~Schuhmacher  \cite{Schu}, familiar with but independent of an early version of this work, has used the machinery of \cite{BKo} to define Hochschild cohomology for morphisms of complex spaces and to give another
proof of the HKR-decomposition theorem in that situation.

One inconvenience of the approach here is the variety of choices involved, from 
the open covering by Stein compact sets, to the local models of the free resolutions, to the glueing data. Although one can track independence of choices step by step, this is indeed cumbersome. Here we use a different approach: to show that our construction is independent of the choices, and thus leads to a well defined object in the derived category, 
we employ {\em categories  of models\/}, in a context as suitable for
our purposes. These are related to but less sophisticated than Quillen's model categories. This technique was already used implicitly in \cite{Fle1} to
deduce that the cotangent complex of a morphism of complex spaces 
is well defined in the derived category. 
Here we formalize the treatment and give complete proofs along
with the application that shows the Hochschild complex 
to be well defined and to behave functorially with respect to morphisms of
complex spaces. 

We use throughout the notations and techniques as set up in
Section 2 of our previous paper \cite{BFl}, including here only quick reviews
of some relevant facts. Moreover, as already mentioned, we treat explicitly just the 
case of morphisms of complex spaces, leaving the simplifications that occur 
in the case of morphisms of schemes to the reader.

\section{The Hochschild complex for complex spaces}

\subsection{Resolvents of complex spaces}

We begin with a short review of the notion of a resolvent of
a morphism of complex spaces, see \cite[p.\ 33]{Fle1},
\cite[2.34]{BFl}, or \cite{Pal}.

\begin{defn}\label{h.1}
A {\em resolvent\/} for a morphism $f:X\to Y$ of complex 
spaces\footnote{A complex space is always viewed as a ringed space $X=(|X|,\calox)$, with
$|X|$ as underlying topological space and $\calox$ as structure sheaf.}, 
or simply for $X$ over $Y$, consists of a triple
$\fX=(X_{*}, W_*,\calr_*)$ given as follows:
\begin{enumerate}
\item
$X_*=(X_\alpha)_{\alpha\in A}$ is the simplicial space
associated to some locally finite covering $(X_i)_{i\in I}$
of $X$ by Stein compact subsets; see \cite[2.2(a)]{BFl}.
In particular, $A$ is the simplicial set of all subsets
$\alpha$ of $I$ with $X_\alpha:=\bigcap_{i\in\alpha}X_i\ne
\emptyset$.

\item
$W_*$ is a {\em smoothing}\/ of $f$, which means that there
is given a factorization
\bdi[small]
X_* && \rInto^{i} && W_*\\
&\rdTo<{f} && \ldTo>{\tilde f}\\
&& Y &&\quad,
\edi
where $W_*=(W_\alpha)_{\alpha\in A}$ is a simplicial system
of Stein compact sets\footnote{Stein compact sets are always assumed to be semianalytic.}
on the same simplicial set as $X_{*}$,
the morphism $i$ is a closed embedding, and the morphism
$\tilde f$ is smooth\footnote{that is, for every point $x\in
W_\alpha$, the analytic algebra $\calo_{W_\alpha,x}$ is
smooth over $\calo_{Y,f(x)}$.}.

\item
$\calr_*$ is a locally free DG algebra%
\footnote{Our DG algebras are always assumed to be concentrated in non-positive degrees.}
over $\calows$, with
each graded component of $\calr_{\alpha}$ a coherent
$\calow_{\alpha}$--module for each simplex $\alpha$, and
there is given a quasiisomorphism $\calr_*\to \caloxs$ of
$\calows$--algebras.
\end{enumerate}

The smoothing $W_*$ is called {\em free} if
\begin{enumerate}
\item[(a)] 
$W_i\subseteq \bbbc^{n_i}\times Y$ for some $n_{i}$, and
$\calo_{W_i}\cong \calo_{\bbbc^{n_i}\times Y}| W_i$;

\item[(b)]  
$W_\alpha\subseteq\prod^Y_{i\in\alpha}W_i$, where $\prod^Y$
denotes the fibre product over $Y$, and
$\calo_{W_\alpha}$ is the topological restriction of the
structure sheaf on $\prod^Y_{i\in\alpha}W_i$ to $W_\alpha$;

\item[(c)] 
the transition maps $W_\beta\to W_\alpha$,
$\alpha\subseteq\beta$, are induced by the corresponding
projections $\prod^Y_{i\in\beta}W_i\to
\prod^Y_{i\in\alpha}W_i$.
\end{enumerate}
The resolvent $\fX$ itself is said to be {\em free} if
$W_*$ is a free smoothing and $\calrs$ is furthermore a {\em
free\/} $\calows$--algebra (see \cite[2.31]{BFl}).
\end{defn}

Each morphism between complex spaces admits such resolvents,
even free ones, with respect to any covering of $X$ by Stein
compact sets; see, for example, \cite[2.11 (a)]{Fle1} or
\cite[2.35]{BFl}. Recall that these resolvents are
constructed inductively along the nerve $A$ of the given
covering. The initial step uses the fact that the ring of
global sections over the Stein compact set $X_{i}$ is
noetherian by Frisch's theorem \cite{Fr}, whence $\calox_{i}$ can be
obtained for some $n_{i}\ge 0$ as a quotient of
$\calo_{{\bbbc^{n_i}\times Y}}$ restricted to a suitable
Stein compact set $W_{i}$ that is smooth over $Y$. The
$\calow_{i}$--algebra resolution $\calr_{i}$ of $\calox_{i}$
can then be constructed \`a la Tate \cite{Tat}; see 
\cite[2.7]{Fle1} and \cite[2.32]{BFl} for the relevant case of DG
algebras over simplicial systems of Stein compact sets. As this result
will be frequently used in the following, we state it along with the
lifting property for free algebra resolutions \cite[2.8]{Fle1}. {\em
By convention, all DG algebras over simplicial systems of Stein
compact sets considered in this paper are assumed to have coherent
homogeneous components concentrated in only non-positive degrees.}
(Recall that the differential {\em increases\/} degree by one, thus, such
a DG algebra is bounded, as a complex, in the direction of the differential.)

\begin{prop}
\label{h6.1}
(1) Given a morphism of DG algebras
$\calas\to\calbs$  on a simplicial system of Stein compact sets there is
a free DG $\calas$-algebra resolution of $\calbs$, that is, a free DG
$\calas$-algebra $\calrs$ together with a quasiisomorphism of DG 
$\calas$-algebras $\calrs\to\calbs$.

(2) Given a commutative diagram in solid arrows 
    \bdi
        \cala_{*}&\rTo & \cala'_{*}\\
        \dTo&\ruDotsto^{\tilde \varphi}&\dTo_{\pi}\\
         \calb_{*}&\rTo^{\vp}& \calb'_{*}
    \edi
of DG algebra morphisms on a
    simplicial scheme of Stein compact sets
    with $ \calb_{*}$ free over $\cala_{*}$ and
    $\pi$ a surjective quasiisomorphism, there exists a
    morphism $\tilde\vp: \calb_{*}\to \cala'_{*}$
    of DG $\cala_{*}$--algebras such that $\pi\circ\tilde\vp
    =\vp$. Moreover, the lifting $\tilde\vp$ of $\vp$
    through $\pi$ is unique in the derived category.
\end{prop}

\begin{proof}
For the proof of (1) we refer the reader to \cite[2.32]{BFl}. To deduce
(2), by induction along the Postnikov tower of the free algebra
$\cala_{*}\to  \calb_{*}$, see \cite[discussion after 2.32]{BFl},
and in view of the structure of free graded modules over a
simplicial DG algebra, \cite[2.13]{BFl}, one may reduce to
the case that $ \calb_{*}$ is obtained from
$\cala_{*}$ by adjunction of a graded free
$\cala_{*}$--module of the form
$p_{\alpha}^{*}{\calf_{\alpha}}$ with $\calf_{\alpha}$
generated in a single degree $k\le 0$ over $\cala_{\alpha}$.
If $e$ is a generator of $\calf_{\alpha}$, then
$\tilde\vp(\partial e)$ is already defined as a section of
$\cala'_{\alpha}$ and $\vp(e)$ is a section in
$ \calb'_{\alpha}$. As $\pi$ is surjective, there exists, by
Theorem B, a section $e'$ in $ \cala'_{\alpha}$ with
$\pi(e')=\vp(e)$. Now $b=\tilde\vp(\partial e)-\partial e'$
is a cycle that maps to zero under $\pi$. As $\pi$ is a
quasiisomorphism, $b=\partial e''$ for some $e''$, and
$\tilde\vp(e)=e'+e''$ yields the desired lifting of
$\vp(e)$.

The last assertion follows as $\tilde\vp = \pi^{-1}\vp$ in 
the derived category.
\end{proof}

\begin{sit}\label{h.2}
Given a resolvent $(X_*,W_*,\calr_*)$ of a morphism of
complex spaces $f:X\to Y$, we can form the simplicial spaces
$X_*\times_Y X_* := (X_\alpha\times_Y X_\alpha)_{\alpha\in
A}$ and $W_*\times_Y W_* := (W_\alpha\times_Y
W_\alpha)_{\alpha\in A}$ of Stein compact sets over $Y$.
Note that $X_*\times_Y X_*$ is the simplicial space
associated to the covering $(X_{i}\times_{Y}X_{i})_{i\in I}$
of a neighbourhood of the diagonal subspace $X\subseteq
X\times_{Y}X$. The simplicial space $W_*\times_Y W_*$
constitutes a smoothing of $X_*\times_Y X_*$ over $Y$. As in
\cite[2.37]{BFl} we set\footnote{another plausible notation,
modeled on the one sometimes used in algebraic geometry, is
$$
\cals_*
=\calr_*{\boxtimes}_{\caloy}
\calr_* := 
p_{1}^{*}\calr_{*}\otimes_{\calo_{W_*\times_{Y} W_*}}
p_{2}^{*}\calr_{*}\,,
$$ 
where $p_{1},p_{2}:W_*\times_{Y} W_*\to W_{*}$ denote the
canonical projections.}
$$
\cals_*:=
\calr_*\tilde\otimes_{\caloy}\calr_* :=
\calr_*\otimes_{\calo_{W_*}}(\calo_{W_*\times_{Y} W_*})
\otimes_{\calo_{W_*}}\calr_*\,,
$$
and note that $\cals_*$ is naturally
a DG algebra over the simplicial space $W_*\times_Y
W_*$. The reader should keep in mind the following facts.
\end{sit}

\begin{lem}
      \label{h.3}
For a point $(w,w')\in W_{\alpha}\times_{Y} W_{\alpha}$, the
complex $\cals_{\alpha, (w,w')}$ is exact if $(w,w')\not\in
X_{\alpha}\times_{Y} X_{\alpha}$, whereas it represents the
derived analytic tensor product $\calo_{X,w}\tilde
{\dotimes}_{\calo_{Y,y}} \calo_{X,w'}$, if $(w,w')\in
X_{\alpha}\times_{Y} X_{\alpha}$ with $y =f(w)=f(w')$. In
particular,

\begin{enumerate}
\item for any simplices $\alpha\subseteq\beta$ and each
$z\in W_\beta\times_YW_\beta$, the transition map
$p_{\alpha\beta}:W_\beta\times_YW_\beta\to
W_\alpha\times_Y W_\alpha$ induces a quasiisomorphism
$\cals_{\alpha, p_{\alpha\beta}(z)} \to \cals_{\beta,z}$;

\item if $X$ is flat over $Y$, then the natural morphism
$\calss\to \calo_{X_{*}\times_Y X_{*}}$ is a
quasiisomorphism.
\end{enumerate}
\end{lem}

\begin{proof}
The given quasiisomorphism $\calr_*\to \calo_{X_*}$ is a
$Y$--flat resolution, whence at each point the DG algebra
$\calss$ represents the derived analytic tensor product
$\calo_{X_*}\tilde{\dotimes}_\caloy\calo_{X_*}$. This
implies the first assertion as well as (1). If $X\to Y$ is
flat, the derived tensor product is represented by the
non-derived one, and the quasiisomorphism $\calr_*\to
\calo_{X_*}$ of $Y$--flat resolutions of $\calo_{X_{*}}$
yields the quasiisomorphism
$$
\calss=\calrs\tilde\otimes_\caloy\calrs\lto
\calo_{X_*}\tilde{\otimes}_\caloy\calo_{X_*}
\cong \calo_{X_{*}\times_Y X_{*}}\,.
$$
\end{proof}

\subsection{The \Cech construction}
We keep the notation from
\ref{h.2}. To construct the Hochschild complex, we use the
\Cech functor as the basic tool to pass from modules on
$X_*$ to modules on $X$. We remind the reader in brief of
the relevant definitions and properties, see also
\cite[2.27-- 2.30]{BFl}.

\begin{sit}
\label{h.4}
Restricting a given $\calox$--module $\calm$ to the Stein
compact sets of the given covering defines the
$\calo_{X_*}$--module $\calm_*=j^{*}\calm$ with
$\calm_\alpha:=\calm|X_\alpha$. This functor is exact and so
induces directly a functor $j^{*}:D(X)\to D(X_*)$ between the
respective derived categories.

To describe a right adjoint, denote by
$j_\alpha:X_\alpha\hto X$ the inclusion and order the
vertices of the covering to associate to any module
$\calm_*$ on $X_*$ the {\em \v{C}ech complex\/}
$\calc^\sbullet(\calm_*)$ with terms
$$
\calc^p (\calm_*):= \prod_{|\alpha|=p}
{j_\alpha}_*(\calm_\alpha)\,,
$$
where the product is taken over all ordered simplices, and
the differential is defined in the usual way by means of the
transition morphisms for $\calm_*$ and the given ordering on
the vertices. The functor $j_*(\calm_*):=
\calh^0(\calc^\sbullet(\calm_*))$ from
$\calo_{X_{*}}$--modules to $\calox$--modules is right
adjoint to $j^{*}$, and the canonical homomorphism of
$\calox$--modules $\calm\to j_*j^{*}(\calm)$ is an
isomorphism. The \Cech functor extends in the usual way to
complexes $\calms$ by taking the total complex of the double
complex $\calc^p(\calm_*^q)$. We note the following
important properties of the \Cech functor:

\begin{enumerate}
    \item 
    \label{C.1}
    $\calc^\sbullet$ is exact, thus can be viewed
    directly as a functor $\calc^\sbullet: D(X_*)\to D(X)$.

    \item
    \label{C.2}
    $\calc^\sbullet$ represents $Rj_*$, the right derived
    functor of $j_*$.
    
    \item
    \label{C.3} 
    The adjunction morphism $\calm\to Rj_*j^{*}\calm
    \cong \calc^\sbullet(\calms)$ is just the \Cech complex of
    sheaves associated to the given covering. It is always a
    quasiisomorphism, that is, the functor $j^{*}:D(X)\to
    D(X_*)$ is fully faithful.
    
    The other adjunction morphism, 
    $j^{*}\calc^\sbullet(\calms)\to \calms$ on $D(X_*)$, is a 
    quasiisomorphism if and only if the transition 
    morphisms $p^{*}_{\alpha\beta}\calm_{\alpha}\to 
    \calm_{\beta}$ are quasiisomorphisms for all pairs of 
    simplices $\alpha\subseteq \beta$.

    \item
    \label{C.4} 
    The terms of the complex $\calc^\sbullet(\calms)$ are
    flat $\calox$--modules whenever $\calm_\alpha$ is flat
    over $\calo_{X_\alpha}$ for each simplex $\alpha$.
\end{enumerate}
\end{sit}

\begin{rem}
We end the general review of the setup recalling the following 
convention. If $\calm,\caln$ are objects of a 
triangulated category $D$, then
$$
\Ext^{n}_{D}(\calm,\caln) := \Hom_{D}(\calm,T^{n}\caln)\quad,\quad n\in\bbbz\,,
$$
with $T$ the translation functor on $D$. 
This convention thus ignores whether $\Ext^{\bdot}$ can be realized
as the cohomology of a ``concrete'' complex, obtained from
resolutions of some kind. This additional feature is however
present in most, if not in all the (unbounded) derived categories 
we work with; see, for example, \cite{Sp}. 
Indeed, these triangulated categories arise 
from abelian categories with enough projectives and/or injectives.
In particular, the sheaves of the form $\cExt^{\bdot}_{X}(\calm,\caln)$
or $\cTor_{\bdot}^{\calox}(\calm,\caln)$ on some space $X$, with 
$\calm,\caln\in D(X)$, exist and are well defined, regardless
of the size of the complexes involved.
\end{rem}

\subsection{The Hochschild complex} 
\label{h.5}
The multiplication map
$\mu:\calss=\calr_*\tilde\otimes_{\caloy}\calr_* \to
\calr_*$ turns $\calrs$ into an $\calss$--algebra, and 
there is a locally free, and by  \ref{h6.1}(1) even a free, DG
$\calss$--algebra resolution
$\calbs$ of
$\calrs$ over $\calss$  that fits into a commutative diagram
\bdi[s=7mm]
&&\calbs \\
&\ruTo  &&\rdTo^{\nu}\\
\calss &&\rTo^{\mu} && \calrs\,,
\edi
where $\nu$ is a quasiisomorphism of DG algebras.

We keep track of these additional data through the following notation. 

\bdfn\label{h.5a}
The quadruple $\fX^{(e)}:=(X_*,
W_*,\calr_*,\calb_*)$ is called an {\em extended resolvent}\/  of
$f:X\to Y$ if $\fX:=(X_*,
W_*,\calr_*)$ is a resolvent and $\calbs$ is a locally free algebra
resolution of $\calrs$ over $\calss=\calr_*\tilde\otimes_{\caloy}\calr_*$. This extended resolvent is
called {\em free}, if $\fX$ is free and $\calbs$ is a free algebra
over $\calss$.
\edfn

By the preceding discussion such an extended free resolvent
always exists. To define Hochschild complexes, free
extended resolvents are sufficient; however to deduce uniqueness we
need the flexibility offered by arbitrary resolvents.

The tensor product of the algebra resolution $\calbs$ over $\calss$ with the composition
$\calss\xto{\mu} \calrs\to\calo_{X_*}$, when restricted topologically to the diagonal
$\Delta:X_{*}\hookrightarrow
X_{*}\times_{Y}X_{*}\hookrightarrow W_{*}\times_{Y} W_{*}$,
leads first to a Hochschild complex of $X_{*}$ over $Y$, and
then through the \v{C}ech functor to a Hochschild complex of
$X$ over $Y$.

\begin{defn}\label{h.6}
We call the DG $\calo_{X_*}$--algebra
$$
\bbbh_{X_*/Y}:= \Delta^{*}(\calb_{*}) =
\left.(\calbs\otimes_\calss \calo_{X_*})\right|_{X_{*}}
$$
a {\em Hochschild complex} of $X_*$ over $Y$. The
associated \v{C}ech complex
$$
\bbbh_{X/Y}:=\calc^\sbullet(\bbbh_{X_*/Y})
$$
of $\calox$--modules will be called a {\em Hochschild
complex} of $X$ over $Y$.
\end{defn}

Unlike the classical Hochschild complex for algebras, the
complexes here are of course not canonically defined as they
depend on the choice of the extended resolvent. However we have
the following result the proof of which will be postponed to
Section 2.

\bthm\label{h.main}
The Hochschild complex is well defined as an object of
the derived category $D(X)$. Moreover, for any commutative diagram
of morphisms of complex spaces
\bdi
X& \rTo^g & X'\\
\dTo && \dTo\\
Y& \rTo & Y'
\edi
there is a functorial morphism
$$
Lg^{*}(\bbbh_{X'/Y'})\to \bbbh_{X/Y}
$$
in $D(X)$.
\ethm

Let us note the following simple properties of Hochschild 
complexes.
\begin{prop}
\label{h.8}    
Let $\bbbh_{X_*/Y}$ and $\bbbh_{X/Y}$ be Hochschild
complexes as just constructed.

\begin{enumerate}
\item
\label{h.8(1)} 
The complex $\bbbh_{X_*/Y}$ is a graded free
$\calox_{*}$--algebra with coherent and locally free
components concentrated in non-positive degrees.

\item
\label{h.8(2)} 
The transition maps $\bbbh_{X_\alpha/Y}|X_\beta\to
\bbbh_{X_\beta/Y}$ are quasiisomorphisms for all simplices
$\alpha\subseteq \beta$. The adjunction morphism
$$
j^{*}\bbbh_{X/Y}\to \bbbh_{X_*/Y}
$$
is a quasiisomorphism.

\item
\label{h.8(3)} 
The cohomology sheaves $\calh^p(\bbbh_{X_*/Y})$ vanish for
$p > 0$ and
$\calh^0(\bbbh_{X_*/Y})\cong \calo_{X_*}$.

\item
\label{h.8(4)}
The complex $\bbbh_{X/Y}$ is locally bounded above with 
coherent cohomology, and its stalks $\bbbh_{X/Y,x}$ are free, so in
particular flat
$\calo_{X,x}$-modules. The cohomology sheaves
$\calh^p(\bbbh_{X/Y})$ vanish for $p> 0$, and
$\calh^0(\bbbh_{X/Y})\cong
\calox$.
\end{enumerate}
\end{prop}

\begin{proof}
The assertions in (\ref{h.8(1)}) are immediate from the
construction. In order to deduce (\ref{h.8(2)}), it suffices
to argue locally. As $\calbs$ is a free algebra over
$\calss$, for each $z\in W_\alpha\times_{Y} W_\alpha$, the
DG $\cals_{\alpha,z}$--module $\calb_{\alpha,z}$ is
projective, resolving $\calo_{X,z}$ via the
quasiisomorphisms $\calb_{\alpha,z} \xto{\simeq}
\calr_{\alpha,z} \xto{\simeq} \calo_{X_{\alpha},z}\cong
\calo_{X,z}$. By \ref{h.2}, for simplices
$\alpha\subseteq\beta$ and every point $z\in
W_\beta\times_{Y} W_\beta$, the corresponding algebra
homomorphism $\cals_{\alpha,p_{\alpha\beta}(z)}\to
\cals_{\beta,z}$ is a quasiisomorphism, where
$p_{\alpha\beta}:W_\beta\times_{Y} W_\beta\to
W_\alpha\times_{Y} W_\alpha$ is as before the transition
map. Applying now \cite[X.66,\S 4 no.3]{Alg}, the morphism
$$
\calb_{\alpha, p_{\alpha\beta}(z)} \to
\calb_{\alpha,p_{\alpha\beta}(z)}
\otimes_{\cals_{\alpha,p_{\alpha\beta}(z)}}\cals_{\beta, z}
$$
is seen to be a quasiisomorphism of
$\cals_{\alpha,p_{\alpha\beta}(z)}$--modules. Hence the
module on the right constitutes along with $\calb_{\beta,z}$
a projective resolution of $\calo_{X,z}$ over $\cals_{\beta,
z}$ and the transition map
$\calb_{\alpha,p_{\alpha\beta}(z)}
\otimes_{\cals_{\alpha,p_{\alpha\beta}(z)}}\cals_{\beta, z}
\to \calb_{\beta,z}$ is consequently a quasiisomorphism of
projective resolutions of $\calo_{X,z}$, thus
$$
\bbbh_{X_\alpha/Y}|X_\beta=\calb_\alpha
\otimes_{\cals_{\alpha}}\calo_{X_\alpha}|
X_\beta \cong p_{\alpha\beta}^{*}(\calb_\alpha)
\otimes_{\cals_\beta}\calo_{X_\beta}|X_\beta \lto \calb_\beta
\otimes_{\cals_\beta}\calo_{X_\beta}|X_\beta=\bbbh_{X_\beta/Y}
$$
is a quasiisomorphism as claimed. The last assertion of
(\ref{h.8(2)}) follows now from \ref{h.4}(\ref{C.3}), and
these arguments prove as well (\ref{h.8(3)}).

It remains to establish (\ref{h.8(4)}). As the given
covering of $X$ is locally finite, $(\bbbh_{X/Y})_{x} =
\calc^\sbullet((\bbbh_{X_*/Y})_{x})$ for every $x\in X$, and
the localized complex is bounded above, whence $\bbbh_{X/Y}$
is locally bounded above. As the stalks
$(\bbbh_{X_*/Y})_{x}$ are free $\calo_{X,x}$--modules by
(\ref{h.8(1)}), the same holds for $(\bbbh_{X/Y})_{x}$. The
quasiisomorphism $ j^{*}\bbbh_{X/Y}\to \bbbh_{X_*/Y} $ 
established in (\ref{h.8(2)}) yields the claims on the 
cohomology sheaves of $\bbbh_{X/Y}$.
\end{proof}

\subsection{The algebra structure on the Hochschild complex}
Note that $\bbbh_{X/Y}$ will no longer necessarily be an
$\calox$--algebra since the \v{C}ech functor is not
compatible with tensor products. However, we will show in
this part that the Hochschild complex is at least an
$\calox$--algebra object in the derived category.

\begin{sit}\label{h.9}
Let $\cala$ be a complex of flat $\calox$-modules. 
Assume given a morphism 
$$
m: \cala\dotimes_\calox\cala=\cala\otimes_\calox\cala\to \cala,
$$
in $D(X)$, called the multiplication. As usual, such a
multiplication is said to be {\em associative} if the
morphisms $m\circ(m\otimes \id_\cala)$ and
$m\circ(\id_\cala\otimes m)$ from
$\cala\otimes_\calox\cala\otimes_\calox\cala$ to $\cala$
are equal. With $\sigma:\cala\otimes_\calox\cala \to
\cala\otimes_\calox\cala$ the Koszul morphism that
interchanges the two factors, the multiplication is (graded)
{\em commutative} if $m=m\circ\sigma$. Finally, a morphism
$\epsilon:\calox\to\cala$ is said to be a (left) {\em
unit}\/ if the morphism $m\circ(\epsilon\otimes \id_\cala)$
from $\calox\otimes_\calox\cala$ to $\cala$ is equal to the
canonical isomorphism $\calox\otimes_\calox\cala\cong
\cala$.

In the following, an object $\cala$ of $D(X)$ consisting of
flat $\calox$-modules and equipped with a commutative and
associative multiplication
$m$ and with a unit $\epsilon$ as above will be called in brief
{\em a commutative $\calox$--algebra} in $D(X)$. Morphisms of such
commutative algebras are introduced in a straightforward
manner. Moreover, if $f:X'\to X$ is a morphism, then $Lf^{*}=f^{*}$
on flat $\calox$-modules, thus it transforms commutative
$\calox$--algebras in
$D(X)$ naturally into commutative $\calo_{X'}$--algebras
in $D(X')$.

In the same way we can introduce commutative
$\caloxs$--algebras in $D(X_*)$. For instance, the
Hochschild complex $\bbbh_{X_*/Y}$ is already equipped with a
commutative, associative and unitary $\calox_{*}$-bilinear
multiplication that is a morphism of complexes, thus it
represents trivially a commutative $\calo_{X_*}$--algebra in
$D(X_*)$.
\end{sit}

We now show that the \v{C}ech functor $\calc^\sbullet$
transforms (commutative) $\calo_{X_{*}}$--algebras in
$D(X_*)$ into $\calox$--algebras in $D(X)$. The 
reader will notice that all we require is indeed the fact 
that the left adjoint $j^{*}$ to $\calc^\sbullet$ commutes 
with (derived) tensor products.

\begin{lem}
\label{h.10}
Let $\calms$, $\calns$ be complexes of 
$\calo_{X_*}$-modules and assume that $\calms$ is
flat over $\caloxs$ and locally bounded above. Then the complex
of $\calox$-modules
$\calcb(\calms)$ is flat,  and there
is a natural morphism
$$
\calcb(\calms)\otimes_\calox\calcb(\calns)\lto
\calcb(\calms\otimes_\caloxs\calns)
$$
in $D(X)$. If, moreover, for all simplices $\alpha \subseteq
\beta$, the transition maps
$p_{\alpha\beta}^{*}\calm_\alpha \to \calm_\beta$
and $p_{\alpha\beta}^{*}\caln_\alpha \to \caln_\beta$ are
quasiisomorphisms, then that natural morphism is an
isomorphism in $D(X)$.
\end{lem}

\begin{proof}
Clearly $\calcb(\calms)$ is flat over $\calox$. 
As $j^{*}$ commutes with tensor products of complexes, one 
obtains a natural morphism
$$
j^{*}\big(\calcb(\calms)\otimes_\calox\calcb(\calns)\big)
\cong j^{*}\calcb(\calms)\otimes_\caloxs
j^{*}\calcb(\calns)\lto \calms\otimes_\caloxs\calns
$$
from the adjunction morphisms $j^{*}\calcb(\calms)\to \calms$
and $j^{*}\calcb(\calns)\to \calns$. Adjunction yields the 
desired natural morphism.

If the transition morphisms for $\calm_{*}$ and $\calns$ are
quasiisomorphisms, then the adjunction morphisms
$j^{*}\calcb(\calms)\to \calms$ and $j^{*}\calcb(\calns)\to
\calns$ are quasiisomorphisms by \ref{h.4}(\ref{C.3}),
whence the displayed morphism is a quasiisomorphism. As
$j^{*}$ is fully faithful on $D(X)$ (see again
\ref{h.4}(\ref{C.3})) the last assertion follows.
\end{proof}

\begin{lem}
\label{h.11}
For each (commutative) flat $\caloxs$--algebra $\calas$ in
$D^-(X_*)$ the associated \Cech complex $\calcb(\calas)$
carries a natural (commutative) $\calox$--algebra structure
in $D^{-}(X)$. If, moreover, the transition morphisms on
$\calas$ are quasiisomorphisms, then the adjunction morphism
$j^{*}\calcb(\calas)\to\calas$ is an isomorphism of algebra
objects in $D^{-}(X_{*})$.
\end{lem}

\begin{proof}
In view of \ref{h.10}, we have a natural morphism
$$
\calcb(\calas)\otimes_{\calox}\calcb(\calas)\lto
\calcb(\calas\otimes_{\calo_{X_{*}}}\calas)
$$
in $D^{-}(X)$. Following this morphism with $\calcb(\mu)$,
where $\mu:\calas\otimes_{\calo_{X_{*}}}\calas\to \calas$
is the given multiplication on $\calas$, defines the
multiplication on $\calcb(\calas)$ in $D^{-}(X)$. That
$\calcb(\calas)$ inherits as well the unit from $\calas$ and 
that these data turn $\calcb(\calas)$ into a (commutative) 
$\calox$--algebra in $D^{-}(X)$ is clear from the naturality 
of the construction, and so is the last assertion.
\end{proof}

Applying the result just established to Hochschild complexes
yields immediately the following.

\begin{prop}\label{h.12}
The complex $\bbbh_{X/Y}$ admits a commutative
$\calox$--algebra structure in the derived category $D(X)$
that does not depend on the choice of $\calbs$. The 
$\calo_{X_{*}}$--algebra $j^{*}\bbbh_{X/Y}$ is canonically 
isomorphic to the $\calo_{X_{*}}$--algebra $\bbbh_{X_{*}/Y}$
in the derived category $D(X_{*})$.\qed
\end{prop}

\begin{rem}\label{h.13}
One can find an explicit morphism of complexes
$$
\bbbh_{X/Y}\times_{\calox} \bbbh_{X/Y} \lto \bbbh_{X/Y}
$$
that represents the multiplication. Indeed, the classical
Alexander Whitney map for simplicial complexes; see, for
example, \cite[VIII 8.5]{MLa}; yields an explicit 
associative unital product
$$
\calc^p(\bbbh_{X_*/Y})\times \calc^q(\bbbh_{X_*/Y}) \lto
\calc^{p+q}(\bbbh_{X_*/Y}\otimes_\caloxs \bbbh_{X_*/Y})
$$
In general, this product is only graded commutative in the
homotopy category of complexes; see loc.cit., VIII 8.7.
Moreover, as we are considering the \Cech complexes of
alternating (or, equivalently, ordered) chains, this
construction depends on an ordering of the vertex set of the
simplices and is not functorial with respect to mappings.
Note, however, that the Alexander-Whitney map becomes
functorial as soon as we replace the alternating \Cech
complex by the larger complex of all \Cech chains as in
\cite{Go}.
\end{rem}

An algebra structure in $D(X)$ induces similar structures in 
cohomology and on various $\Tor$ or $\Ext$ groups or sheaves.
More precisely, we have the following standard application.

\begin{lem}\label{h.14}
An $\calox$--algebra structure on a complex $\cala\in D^-(X)$ of
flat
$\calox$-modules induces a natural commutative graded
$H^0(X,\calox)$--algebra structure on $H^{\bdot}(X,\cala)$.
Moreover, for every complex $\calm\in D(X)$, the groups
$$
\Tor_{\bdot}^X(\calm, \cala)\quad\mbox{and}\quad
\Ext^{\bdot}_X(\cala, \calm)
$$
carry natural graded module structures over $H^{\bdot}(X,\cala)$.
Similarly, the sheaf $\calh^{\bdot}(\cala)$ is a graded
commutative $\calox$--algebra, and the sheaves
$$
\cTor_{\bdot}^X(\calm, \cala)\quad\mbox{and} \quad\cExt^{\bdot}_X(\cala, \calm)
$$
are modules over it.
\end{lem}

\begin{proof}
Note that $H^{\bdot}(X,\cala) = \Ext^{\bdot}_{X}(\calox, \cala)$ by 
definition. If $f,g:\calox\to \cala$ are morphisms in 
$D^{-}(X)$ representing elements of this graded group of 
degree $|f|,|g|$ respectively, then
the composition
$$
\calox\xto{\cong} \calox\otimes_{\calox}\calox
\xto{f\dotimes g}\cala\otimes_{\calox}\cala\xto{\mu}\cala
$$
defines the product $fg\in H^{|f|+|g|}(X,\cala)$. 
It is easy to see that this multiplication is associative,
with unit the image of $1$ under
$H^{\bdot}(X,\calox)\xto{\epsilon} H^{\bdot}(X,\cala)$. To show that it
is graded commutative, note that, with $T$ the translation funtor on
$D(X)$, the diagram
\bdi[s=8mm]
T^n \cala\otimes_\calox T^m\cala &\rTo^{\cong} & T^{n+m}
(\cala\otimes_\calox\cala)\\
\dTo^\sigma_{\cong} &(-1)^{mn}& \dTo^\sigma_{\cong}\\
T^m\cala\otimes_\calox T^n\cala &\rTo^{\cong} & T^{n+m}
(\cala\otimes_\calox\cala)
\edi
commutes up to the sign $(-1)^{nm}$, where $\sigma$ denotes
as before the Koszul map interchanging the two factors. Taking
cohomology, it follows easily that $H^{\bdot}(X,\cala)$ is graded
commutative. The remaining assertions are left to the reader
as an exercise.
\end{proof}

\section{Hochschild (co-)homology of complex spaces} 
\subsection{Definition and basic properties}
With Hochschild
complexes at our disposal, it is now immediate how to define
Hochschild (co-)homology.

\begin{defn}\label{h.15}
Let $\calm\in D(X)$ be a complex of $\calox$--modules. The groups
\begin{align*}
    \HH_{\bdot}^{X/Y}(\calm) &:=\Tor_{\bdot}^X(\bbbh_{X/Y}, \calm) :=
    H^{-\bdot}(X,\bbbh_{X/Y}\dotimes_{\calox} 
    \calm)\quad\text{and}\\    
    \HH^{\bdot}_{X/Y}(\calm) &:=\Ext_X^{\bdot}(\bbbh_{X/Y}, \calm)
\end{align*}
are called the {\em Hochschild homology}, resp.\ {\em Hochschild
cohomology} of $\calm$. Similarly we introduce {\em Hochschild
(co-)homology sheaves\/}
$$
\cHH_{\bdot}^{X/Y}(\calm):=\cTor_{\bdot}^X(\bbbh_{X/Y}, \calm)
\quad\mbox{and}\quad
\cHH^{\bdot}_{X/Y}(\calm):=\cExt_X^{\bdot}(\bbbh_{X/Y}, \calm)\,.
$$

If $Y$ is just a point, we write simply $\HH^{X}_{\bdot}(\calm), \HH_{X}^{\bdot}(\calm)$,
and so forth, and call these groups or sheaves {\em absolute Hochschild
(co-)homology\/} of $X$.
\end{defn}

The following properties are immediate from the definition and the
previous results.

\begin{prop}\label{h.16} With Hochschild (co-)homology as 
just introduced, one has
\begin{enumerate}
\item $\HH^{\bdot}_{X/Y}$ is a cohomological functor, that is,
every short exact sequence\footnote{The reader may, of course, substitute ``distinguished triangle'' for ``short exact sequence'', if needed.} of complexes of $\calox$--modules
$0\to\calm'\to\calm\to\calm''\to 0$ induces a long exact
cohomology sequence
$$
\cdots \to \HH^{p-1}_{X/Y}(\calm'')\to \HH^{p}_{X/Y}(\calm')\to
\HH^{p}_{X/Y}(\calm)\to \HH^{p}_{X/Y}(\calm'')\to\cdots\,,
$$
and these long exact sequences depend functorially on the
given short exact sequences. Analogously, $\cHH^{\bdot}_{X/Y}(-)$ is a
cohomological functor, and similarly
$\HH_{\bdot}^{X/Y}(-)$ and $\cHH_{\bdot}^{X/Y}(-)$ are homological functors on
$D(X)$.

\item If $\calm\in D_{coh}^-(X)$, then the sheaves
$\cHH_p^{X/Y}(\calm)$ are coherent
for each $p\in\bbbz$. 
Similarly, the sheaf $\cHH^p_{X/Y}(\calm)$ is coherent for 
$\calm\in D_{coh}^+(X)$ and any $p$.

\item For every $\calox$--module $\calm$ and any $p<0$, the objects
$$
\HH^p_{X/Y}(\calm)\,,\quad
\cHH^p_{X/Y}(\calm)\,,
\quad\text{and}\quad \cHH_p^{X/Y}(\calm)\,,
$$
vanish. Moreover, $\HH^0_{X/Y}(\calm)\cong \Hom_X(\calox,\calm)=H^{0}(X,\calm)$, and
$$
\cHH_0^{X/Y}(\calm)\cong \calm\cong \cHH^0_{X/Y}(\calm)
$$

\item 
The Hochschild homology $\HH_{\bdot}^{X/Y}(\calox)$ carries a
natural graded commutative $H^0(X,\calox)$--algebra
structure, and for every complex $\calm\in D(X)$ the groups
$\HH_{\bdot}^{X/Y}(\calm)$ and $\HH^{-\bdot}_{X/Y}(\calm)$ are graded
modules over this algebra. Similarly, $\cHH_{\bdot}^{X/Y}(\calox)$
is a graded commutative sheaf of $\calox$--algebras, and
$\cHH_{\bdot}^{X/Y}(\calm)$, $\cHH^{-\bdot}_{X/Y}(\calm)$ are graded
modules over it.
\end{enumerate}
\end{prop}

\begin{proof}
Taking into account \ref{h.8} and \ref{h.12},
(1) through (3) are standard properties of $\Ext$ and $\Tor$--groups or sheaves.
Finally, \ref{h.12} and \ref{h.14} imply (4).
\end{proof}

\brem
Note that the Hochschild homology groups
$\HH_{p}^{X/Y}(\calox)$ will generally not be zero for $p<0$, whence the graded
algebra structure may involve both positive and negative degrees.
If $\calm$ is an $\calox$--module with proper support of dimension $d$, then
$\HH_{p}^{X/Y}(\calm)=0$ for $p<-d$, while $\HH_{-d}^{X/Y}(\calm)\cong
H^{d}(X,\calm)$. This follows from the hyper(co-)homology spectral sequence
\begin{align*}
E^{2}_{p,q}=H^{-p}(X,\cHH_{q}^{X/Y}(\calm))\Longrightarrow \HH_{p}^{X/Y}(\calm)
\end{align*}
and \ref{h.8}(\ref{h.8(3)}).
\erem

\subsection{Commutativity of Hochschild (co-)homology} 
We will now show that the Hochschild cohomology
$\HH^{\bdot}_{X/Y}(\calox)$ also admits a natural graded
commutative algebra structure, induced from the Yoneda 
or composition product. For arbitrary (affine) algebras $A\to B$, 
this is Gerstenhaber's famous theorem \cite{Ger}. 
The proof here applies the inspired treatment of the Eckmann-Hilton
argument in \cite{SA}.
The multiplicative structure in question is based on the
following isomorphisms.

\begin{lem}\label{h.17}
(1) There are canonical isomorphisms
$$
\HH^{\bdot}_{X/Y}(\calox)\cong \Ext^{\bdot}_\calss(\calrs,\calrs) 
\cong\Ext^{\bdot}_\calss(\calbs,\calbs)
$$
and the Yoneda product on the last two terms endows thus 
$\HH^{\bdot}_{X/Y}(\calox)$
with the (same) structure of a graded algebra.

(2) For every complex of $\calox$--modules $\calm$, there are
natural isomorphisms
$$
\HH^{\bdot}_{X/Y}(\calm)\cong \Ext^{\bdot}_\calss(\calrs,\calms) 
\cong \Ext^{\bdot}_\calss(\calbs,\calms)
$$
and the action of the corresponding Yoneda $\Ext$--algebra through 
the contravariant argument of  $\Ext^{\bdot}_\calss(?,\calms)$ realizes 
$\HH^{\bdot}_{X/Y}(\calm)$ as a graded
right module over $\HH^{\bdot}_{X/Y}(\calox)$.

(3)  With $\calm$  as in (2), there are natural isomorphisms
$$
\HH_{\bdot}^{X/Y}(\calm)\cong
\Ext^{-\bdot}_\caloxs(\caloxs,\calrs\dotimes_\calss\calms)\cong
\Ext^{-\bdot}_\calrs(\calrs,\calrs\dotimes_\calss\calms)\,,
$$
and the pairing
$$
\Ext^{\star}_\calss(\calrs, \calrs) \times 
\Ext^{\bdot}_\calrs(\calrs, \calrs\dotimes_\calss\calms) 
\lto
\Ext^{\star+\bdot}_\calrs(\calrs, \calrs\dotimes_\calss\calms)
$$
that sends $(\varphi, \psi)$ to
$(\varphi\dotimes_{\calss} \id_\calms)\circ\psi$
endows $\HH_{-\bdot}^{X/Y}(\calm)$ with the structure of a graded left
$\HH^{\bdot}_{X/Y}(\calox)$--module.
\end{lem}

\begin{proof}
To establish the isomorphism in (2), note first that for
all simplices $\alpha\subseteq \beta$ the transition maps
$\bbbh_{X_\alpha/Y}|X_\beta \to \bbbh_{X_\beta/Y}$ are
quasiisomorphisms by \ref{h.8}(\ref{h.8(2)}) and so, applying
\cite[2.30 (2)]{BFl}, there is an isomorphism
$$
\HH^{\bdot}_{X/Y}(\calm)= \Ext^{\bdot}_X(\bbbh_{X/Y}, \calm)\cong
\Ext^{\bdot}_{X_*}(\bbbh_{X_*/Y}, \calms).
$$
According to \cite[2.25 (2)]{BFl}, the term on the right is isomorphic
to $\Ext^{\bdot}_\calss(\calbs,\calms)$. As $\calbs\to\calrs$ is a
quasiisomorphism, the latter group is as well isomorphic to
$\Ext^{\bdot}_\calss(\calrs,\calms)$.
The isomorphism in (1) is just the special case $\calm =\calox$,
and in terms of the $\Ext$--groups, the assertions on the multiplicative 
and linear structures hold for any triangulated category.

Finally, (3) is deduced in the same fashion, and we leave the details to the reader.
\end{proof}
\begin{rem}
We remind the reader that a DG algebra is usually {\em not projective\/} as a DG module over itself or a subalgebra; see \cite[Ex. 2.20]{BFl}. However, the algebras in question here admit projective approximations according to (loc.cit). So, if $\calps\to\calbs$ is a projective approximation of the DG $\calss$--module $\calbs$ and $\calqs\to\calrs$ one of $\calrs$ as a DG module over itself, then we may realize the groups in the preceding lemma as (co-)homology groups of complexes:
\begin{align*}
\HH^{\bdot}_{X/Y}(\calox)&\cong H^{\bdot}(\Hom_{\calss}(\calps, \calps))\cong H^{\bdot}(\Hom_{\calss}(\calps, \caloxs))\\
\HH^{\bdot}_{X/Y}(\calm)&\cong H^{\bdot}(\Hom_{\calss}(\calps, {\widetilde\calm}_{*}))\\
\HH_{\bdot}^{X/Y}(\calm)&\cong H^{-\bdot}(\Hom_{\calrs}(\calqs,\calps\otimes_{\calss} {\widetilde\calm}_{*}))
\end{align*}
where ${\widetilde\calm}_{*}$ is a $W_{*}$-acyclic resolution of $\calms$; see \cite[2.23]{BFl}.

In particular, the pairing in \ref{h.17}(3) can be expressed in terms of complexes 
as the pairing
\begin{align*}
\bdi
H^{\star}(\Hom_{\calss}(\calps,\calps))\times H^{\bdot}(\Hom_{\calrs}(\calqs, \calps\otimes_{\calss}{\widetilde\calm}_{*}))\\
\dTo\\
H^{\star+\bdot}(\Hom_{\calrs}(\calqs, \calps\otimes_{\calss}{\widetilde\calm}_{*}))
\edi
\end{align*}
that sends $(\vp,\psi)$ to $(\vp\otimes_{\calss}\id_{{\widetilde\calm}_{*}})\circ\psi$.
In this explicit description, it becomes clear, for example, that
$\HH^{\bdot}_{X/Y}(\calox)\cong \Ext^{\bdot}_{\calss}(\calrs, \calrs)$ acts naturally on
$\HH^{X/Y}_{-\bdot}(\calm)\cong \Ext^{\bdot}_{\calrs}(\calrs,\calrs\dotimes_{\calss} \calm)$ only through the covariant argument. 
\end{rem}
\begin{cor} 
\label{h.17.5}
For any morphism $f:X\to Y$ and any complex $\calm$ 
of $\calox$--modules, there are natural maps
\begin{align*}
&\HH_{\bdot}^{X/Y}(\calm)\lto \Tor_{\bdot}^{\calo_{X\times_{Y}X}}(\calox,\calm)\\
&\Ext^{\bdot}_{X\times_{Y}X}(\calox,\calm)\lto \HH^{\bdot}_{X/Y}(\calm)\,.
\end{align*}
For $\calm=\calox$, the second map becomes a homorphism of graded algebras,
and the given maps are homomorphisms of graded modules over it.
\end{cor}

\begin{proof}
We show the existence of the second map, leaving the analogous argument for the first one to the reader. From \ref{h.17}(2) we have first $\HH^{\bdot}_{X/Y}(\calm)\cong  \Ext^{\bdot}_{\calss}(\calbs, \calms)$, and then $ \Ext^{\bdot}_{\calss}(\calbs, \calms)\cong
\Ext^{\bdot}_{\calss}(\caloxs, \calms)$ as $\calbs\to\caloxs$ is a quasiisomorphism.
The algebra morphism $\calss\to \caloxs\tilde\otimes_{\caloy}\caloxs$ provides a 
``forgetful'' algebra homomorphism 
$$
\Ext^{\bdot}_{ \caloxs\tilde\otimes_{\caloy}\caloxs}(\caloxs, \calms)\to \Ext^{\bdot}_{\calss}(\caloxs, \calms)\,,
$$
and finally $\Ext^{\bdot}_{ \caloxs\tilde\otimes_{\caloy}\caloxs}(\caloxs, \calms)\cong
\Ext^{\bdot}_{X\times_{Y}X}(\calox,\calm)$, in view of \ref{h.4}(\ref{C.3}).
\end{proof}

\bthm\label{h.18} 
The Hochschild cohomology $\HH^{\bdot}_{X/Y}(\calox)$ is a graded commutative 
$\Gamma(X,\calox)$--algebra with respect to the Yoneda
product.
\ethm

\begin{proof}
To show commutativity of the Yoneda product, 
we apply \cite[Theorem 1.7]{SA} to the derived category
$\catc=D^-(\calss)$ of DG $\calss$-modules. To this end, we
show first that $D^-(\calss)$ is a suspended monoidal category in the sense of
Definition 1.4 in \cite{SA}. In fact,
$\calss$ has a $\calrs$-algebra structure from the left and
from the right, so the analytic tensor product
$\calss\tilde\otimes_\calrs\calss$ is an $\calss$-bimodule. Given
$\calss$-modules $\calms$ and $\calns$,
their analytic tensor product
$$
\calms\tilde\dotimes\calns
:=\calms\dotimes_\calss(\calss\tilde\otimes_\calrs\calss)
\dotimes_\calss\calns
$$
carries a natural $\calrs$-structure from the right and from
the left and thus it admits again an $\calss$-structure. As
$\calss\tilde\otimes_\calrs\calss$ is a flat $\calss$-module
from the right we have, for $\calns=\calrs$,
$$
(\calss\tilde\otimes_\calrs\calss)\dotimes_\calss \calrs=
(\calss\tilde\otimes_\calrs\calss)\otimes_\calss \calrs
\cong \calss\tilde\otimes_\calrs  \calrs\cong \calss.
$$
As a consequence, there are isomorphisms
$$
\varrho_\calms: \calms\tilde\dotimes\calrs\to\calms
\quad\mbox{and, similarly, }\quad
\lambda_\calms: \calrs\tilde\dotimes \calms\to\calms
$$
that are functorial in $\calms\in D^-(\calss)$.   
Hence $\calrs$ is a ``unit'' in the category
$D^-(\calss)$ with respect to $\tilde\dotimes$. 
As $\tilde\dotimes$ satisfies further the usual associativity and
commutativity properties of tensor products, it follows that
$D^-(\calss)$ is indeed a suspended monoidal category. 
Applying the main result of \cite{SA}, the commutativity of 
$$
\Ext^{\bdot}_\calss(\calrs,\calrs)=
\bigoplus_{i\in\bbbz}\Hom_{D^-(\calss)}(\calrs, T^i\calrs)
$$
follows; here $T$ denotes again the translation functor on the derived
category. For the convenience of the reader we give a short
outline of the argument: it is immediately
seen from the construction that $\varrho_\calrs=\lambda_\calrs$.
Given morphisms $f:\calrs\to \calrs$, of degree $|f|$, and $g:\calrs\to
\calrs$, of degree $|g|$, in the derived category, functoriality of
$\varrho_\calms$ and $\lambda_\calns$ shows that the diagram
\bdi
\calrs\tilde\dotimes \calrs &\rTo^{1\otimes f}
&\calrs\tilde\dotimes\calrs & \rTo<{g\otimes 1} &
\calrs\tilde\dotimes \calrs\\
\dTo<{\lambda_\calrs} && 
\dTo<{\lambda_\calrs}>{=\varrho_\calrs}
&&\dTo>{\varrho_\calrs}\\
\calrs&\rTo^f & \calrs &\rTo^{g} & \calrs\\
\edi
commutes in $D^-(\calss)$. In other words, identifying $\calrs$
with $\calrs\tilde\dotimes\calrs$ via $\varrho_\calrs=\lambda_\calrs$,
the morphisms $g\circ f$ and $g\otimes f=(g\otimes 1)\circ (1\otimes
f)$ of degree $|f|+|g|$ are equal in $D^-(X)$. Similarly,
$f\circ g$ is equal to
$(1\otimes f)\circ (g\otimes 1) = (-1)^{|f||g|}g\otimes f$, whence
graded commutativity follows. 
\end{proof}

Given a composable pair of morphisms $f:X\to Y$ and $g:Y\to Z$,
of complex spaces, there is a natural morphism 
$\bbbh_{X/Z} \to \bbbh_{X/Y}$ in the derived
category of $X$; see \ref{h.main}. It induces a natural map
$\HH_{X/Y}^{\bdot}(\calm)\to\HH_{X/Z}^{\bdot}(\calm)$ for every complex 
$\calm$ of $\calox$--modules. 

We postpone the proof of the following result to \ref{m.11}.

\begin{prop}\label{h.23}
The natural map $\HH_{X/Y}^{\bdot}(\calox)\to\HH_{X/Z}^{\bdot}(\calox)$ is
a homomorphism of graded commutative algebras, and the induced
map $\HH_{X/Y}^{\bdot}(\calm)\to\HH_{X/Z}^{\bdot}(\calm)$ is a homomorphism of modules over this algebra homomorphism.

In particular, there is a homomorphism of graded commutative algebras from the relative Hochschild cohomology of $X$ over $Y$ to the absolute Hochschild cohomology,  $\HH_{X/Y}^{\bdot}(\calox)\to\HH_X^{\bdot}(\calox)$.
\end{prop}

\begin{rem}
We end this section raising a subtle point that already appears in the affine case 
and provides further evidence that the derived version of Hochschild (co-)homology as originally proposed by Quillen in \cite{Qui3} for affine algebras and employed here for complex spaces is indeed appropriate. 

For an associative, not necessarily commutative, algebra $A\to B$ over the 
commutative ring $A$, let $B^{e}:= B^{op}\otimes_{A}B$ be the enveloping algebra 
and $\bbbb_{\bdot} :=  \bbbb_{\bdot} (B/A)$ the bar resolution of $B$ over $A$, 
a resolution of $B$ as a (right) module over $B^{e}$ via the multiplication map. 
The cohomology 
$H^{\bdot}(\Hom_{B^{e}}(\bbbb_{\bdot},B))\cong H^{\bdot}(\Hom_{B^{e}}(\bbbb_{\bdot}, \bbbb_{\bdot}))$ constitutes then the classical, non-derived Hochschild cohomology of $B$ over $A$ and the composition of endomorphisms of $\bbbb_{\bdot}$ induces the Yoneda product on it. Remarkably, this product on classical Hochschild cohomology is always graded commutative, a fact first observed by Gerstenhaber \cite{Ger}, who also supplied a simple direct proof in \cite{Ger3}.

If $\bbbp_{\bdot}\to B$ is a projective resolution of $B$ as a $B^{e}$--module, then there exists a comparison map $\bbbp_{\bdot}\to \bbbb_{\bdot}$ of complexes of $B^{e}$--modules over the identity on $B$. This comparison map is, of course, unique in the homotopy category and induces thus a natural homomorphism of graded $A$--algebras 
\begin{align*}
\alpha: H^{\bdot}(\Hom_{B^{e}}(\bbbb_{\bdot},B)) \lto \Ext^{\bdot}_{B^{e}}(B,B)
\end{align*}
that is an isomorphism as soon as $B$ is projective over $A$, as then $\bbbb_{\bdot}$ itself constitutes already a projective $B^{e}$--module resolution of $B$.

On the other hand, consider Quillen's \cite{Qui3} derived version of Hochschild cohomology, constructed by choosing a free (associative) DG $A$--algebra resolution $R$ of $B$ over $A$, and then
taking cohomology, 
$$
\Ext^{\bdot}_{B\dotimes_{A}B}(B,B) := H^{\bdot}(\Hom_{R^{op}\otimes_{A}R}(R,B))\,.
$$
There is as well a natural comparison map
\begin{align*}
\beta: \Ext^{\bdot}_{B^{e}}(B,B)\lto \Ext^{\bdot}_{B\dotimes_{A}B}(B,B)
\end{align*}
that is an $A$--algebra homomorphism with respect to the Yoneda product on either side.
This homomorphism is an isomorphism, by \cite{Qui3}, as soon as the ``transversality condition''
\begin{align}
\label{ortho}
\tag{$\dagger$}
\Tor^{A}_{i}(B,B) = 0\quad\text {for $i>0$}
\end{align}
is satisfied, in particular, when $B$ is flat over $A$. The Eckmann-Hilton argument, as formulated by Suarez-Alvarez \cite{SA}, that we employed above yields graded commutativity of $\Ext^{\bdot}_{B\dotimes_{A}B}(B,B)$ for any algebra $B$ over $A$.

However, the argument does not establish commutativity of $ \Ext^{\bdot}_{B^{e}}(B,B)$ in general! The rather delicate point is that, without transversality as in (\ref{ortho}), the object $B$ is {\em not necessarily a unit\/} for the derived bifunctor 
$$
?\dotimes_{B}?:D(\Mod B^{e})\times D(\Mod B^{e})\lto D(\Mod B^{e})\,.
$$
Indeed, the underlying bifunctor $?\otimes_{B}?$ on the module category commutes with all colimits, that is, it is right exact and respects all direct sums in either argument. The left derived tensor product 
$B\dotimes_{B}B$ is consequently always represented by the total complex of $\bbbp_{\bdot}\otimes_{B}\bbbp_{\bdot}$, where $\bbbp_{\bdot}$ is a projective  $B^{e}$--module resolution of $B$ as above. While this object comes equipped with the natural augmentation 
$$
B\dotimes_{B}B=\bbbp_{\bdot}\otimes_{B}\bbbp_{\bdot}\lto H_{0}(\bbbp_{\bdot}\otimes_{B}\bbbp_{\bdot})\cong B\,,
$$
which equals $\rho_{B}=\lambda_{B}$ in the notation adapted from the proof of \ref{h.18} or \cite{SA}, this augmentation need not be a quasiisomorphism, as, say, the case of the homomorphism of commutative rings $A=K[y]\to B=K[x,y]/(x^{2},xy)$, $K$ some field, already demonstrates. For a different view of the importance of some kind of transversality conditions such as (\ref{ortho}) in this context, see also \cite{Schw}.
\end{rem}

\subsection{The case of a flat morphism}
Next we will show that for a flat morphism
Hochschild (co-)homology is nothing but the usual (Ext-)Tor- algebra
of the diagonal in $X\times_YX$. The reader may wish to compare this to the corresponding result for quasi-projective schemes over a field, as described in \cite{Sw}.

\begin{prop}\label{h.19}
For any  flat morphism $X\to Y$ of complex spaces, there is a canonical isomorphism of
$\calox$--algebras in $D(X)$
\begin{equation}
  \label{1.13.1}
  \bbbh_{X/Y} \cong \calox{\dotimes}_{\calo_{X\times_YX}} \calox\cong (L\Delta^{*})\Delta_{*}\calox,
\end{equation}
where we consider $\calox$ as a module on $X\times_YX$ via
the diagonal embedding $\Delta:X\hto X\times_YX$.
Accordingly, for every complex of $\calox$--modules $\calm$, the comparison 
maps from \ref{h.17.5} become isomorphisms,
\begin{equation}\label{1.13.2}
\begin{array}{c}
\HH^{\bdot}_{X/Y}(\calm)\cong \Ext^{\bdot}_{X\times_YX}(\calox, \calm)\\[3pt]
\HH_{\bdot}^{X/Y}(\calm)\cong \Tor_{\bdot}^{X\times_YX}(\calox, \calm).
\end{array}
\end{equation}
For $\calm=\calox$, these isomorphisms are compatible with the algebra
structures on either side.
\end{prop}

\begin{proof}
First note that by definition
$$
\bbbh_{X_*/Y}=\calbs\otimes_\calss\calo_{X_*}
\cong
(\calbs\otimes_\calss\calo_{X_*\times_YX_*})
\otimes_{\calo_{X_*\times_YX_*}} \calo_{X_*}.\leqno (*)
$$
As $X\to Y$ is flat, the canonical map $\calss\lto\calo_{X_*\times_YX_*}$
is a quasiisomorphism, by \ref{h.3} (2). Tensoring with $\calbs$ from the left we
get that
$$
\calbs\to \bar\calbs:=\calbs\otimes_\calss\calo_{X_*\times_YX_*}
$$
is as well a  quasiisomorphism (cf.\ \cite[2.25 (3)]{BFl}). This shows
that $\bar\calbs$ is a flat resolution of $\calo_{X_*}$ over 
$\calo_{X_*\times_YX_*}$ and so the right hand side of $(*)$ represents
$\calo_{X_*}\dotimes_{\calo_{X_*\times_YX_*}}\calo_{X_*}$. Taking the
\v{C}ech complex, (\ref{1.13.1}) follows.

The remaining isomorphisms (\ref{1.13.2}) and the remark concerning algebra structures are immediate consequences
of (\ref{1.13.1}), the definitions and \ref{h.17.5}. 
\end{proof}

\brem
\label{h.19a}
Observe that \ref{h.19} applies in particular to the absolute case,
when $Y$ is a simple point or $X$ a scheme over a field. Thus, if one were 
interested only in the absolute case, one could {\em define\/} Hochschild 
(co\nobreak-)homology through these $\Tor$-- or $\Ext$--groups, the point of view 
taken, say, in \cite{Sw}. 

However, the definition given here allows to define (absolute) 
Hochschild {(co\nobreak-)ho\-m\-ology} as well for arbitrary schemes (over $\Spec\bbbz$):
while the structure morphism of such a scheme to $\Spec\bbbz$ generally 
need not be flat, the construction via resolvents still applies
mutatis mutandis, replacing Stein compacts by affine schemes. 
The comparison maps as in \ref{h.17.5} still exist, but they will not be isomorphisms
in general.
\erem

\section{Hochschild cohomology and the centre of the derived category}
The aim of this section is to show that there exists a characteristic homomorphism
of graded commutative algebras from Hochschild cohomology of $f:X\to Y$ to the
graded centre of the derived category of $X$. We begin with a review of some terms.

\subsection{Centre of a category}
Let $\catc$ be any category. Recall that a natural transformation
$f:\id_\catc\to\id_\catc$, or endomorphism of the identity functor, 
consists of morphisms $f_M:M\to M$, one for each $M$  in $\catc$,
such that for every morphism
$\alpha:M\to N$ in $\catc$ the diagram
$$
\bdi[h=7mm]
M &\rTo^{\alpha} &N\\
\dTo<{f_M} && \dTo_{f_N}\\
M &\rTo^{\alpha} &N
\edi\leqno (*)
$$
commutes. If $\catc$ is {\em small\/}, these endomorphisms of $\id_{\catc}$ form a set, 
the so-called {\em centre}\/
$\fZ(\catc)$ of $\catc$,
$$
\fZ(\catc):=\Hom (\id_\catc,\id_\catc),
$$
see e.g.\ \cite{MLa}.
Composition of  morphisms of functors provides a product on
$\fZ(\catc)$ with the identity transformation as unit. If $f,g$
are elements of $\fZ(\catc)$ then applying
$(*)$ to the case $M=N$ and $\alpha=g_M$ yields that
the centre is always commutative with respect to this product.

If $\catc$ is furthermore a {\em $K$--linear\footnote{that is, 
$K$ is a commutative ring, each $\Hom$-set is a $K$--module, 
and composition is $K$--bilinear.}\/} category, then endomorphisms 
of the identity functor form themselves a $K$--module and the 
centre $\fZ(\catc)$ comes equipped with the structure of a 
commutative $K$--algebra.

\subsection{Graded centre of a triangulated category}
Now let $\catc = (\catc, T,\Delta)$ be a triangulated category with translation
functor $T$ and collection of distinguished triangles $\Delta$. 
The category $\catc$ is in particular graded by $T$, and we can 
consider more generally the abelian groups 
$\fZ_{gr}^n:=\Hom^T_\catc(\id_\catc, T^n)$ of all natural
transformations, or morphisms of functors, $f:\id_\catc\to T^n$ 
that {\em anticommute with the shift functor}, so that for
every object $M$ of $\catc$ we have
$$
f_{T^pM}= (-1)^{pn}T^p(f_{M}): T^pM\lto T^{p+n}M,
$$
or, in brief, $fT^p=(1)^{pn}T^pf$. The direct sum
$$
\fZ_{gr}^{\bdot}(\catc)=\bigoplus_{n\in\bbbz}\Hom^T_\catc(\id_\catc, T^n),
$$
is the {\em graded centre} of $\catc$.
It is a graded commutative ring: the product of two elements $f\in
\Hom^T_\catc(\id_\catc, T^n)$ and $g\in \Hom^T_\catc(\id_\catc, T^m)$
is given by $T^m(f)\circ  g\in \Hom^T_\catc(\id_\catc, T^{n+m})$.

What precisely is an element of $\fZ_{gr}^{n}(\catc)$? This question is answered 
by the following explicit description.

\begin{lemma}
An element  $f\in \fZ_{gr}^{n}(\catc)$ is represented by a collection of morphisms $f_{\calm}:\calm\to T^{n}\calm$, equivalently, elements  $f_{\calm}\in \Ext^{n}_{\catc}(\calm,\calm)$, one for each object $\calm$ in $\catc$, such that for each distinguished triangle
\begin{align*}
\calm'\xto{u}\calm\xto{v}\calm''\xto{w} T\calm\,,
\end{align*}
thus, $(u,v,w)$ in $\Delta$, the following diagram commutes:
\begin{align*}
\bdi[h=7mm]\tag{$**$}
\calm' &\rTo^{u} &\calm&\rTo^{v}&\calm''&\rTo^{w}&T\calm'\\
\dTo^{f_{\calm'}} && \dTo^{f_{\calm}}&& \dTo_{f_{\calm''}}&& \dTo_{T(f_{\calm'})}\\
T^{n}\calm' &\rTo^{T^{n}u} &T^{n}\calm&\rTo^{T^{n}v}&T^{n}\calm''
&\rTo^{(-1)^{n}T^{n}w}&T^{n+1}\calm'
\edi
\end{align*}
\end{lemma}

\begin{proof}
Any morphism $u:\calm'\to \calm$ in $\catc$ occurs as the first component in a distinguished triangle by one of the axioms of  a triangulated category, whence commutation of the leftmost square subsumes that $f$ is a natural transformation from $\id_{\catc}$ to $T^{n}$.

The bottom row in the diagram $(**)$ is again a distinguished triangle in $\catc$ by the translation axiom for a triangulated category, and anticommutation of $f$ with the translation functor just means that the rightmost square commutes as well. Thus, the elements of  $\fZ_{gr}^{n}(\catc)$ are precisely those morphisms from $\id_{\catc}$ to $T^{n}$ that are natural with respect to distinguished triangles. 
\end{proof}

\begin{rem}
Note that applying $\Hom_{\catc}(\caln,?)$ or 
 $\Hom_{\catc}(?,\caln)$ to the diagram $(**)$ returns the standard long exact sequences for the $\Ext$-groups. Invoking Yoneda's Lemma, commutativity of $(**)$ now simply means that the family $(f_{\calm})_{\calm}$ is functorial with respect to long exact sequences of $\Ext$--groups.
\end{rem}

\begin{sit}
Given an object $\calm$ in $\catc$, the {\em evaluation map\/},
\begin{align*}
ev_{\calm}:  \fZ_{gr}^{\bdot}(\catc)\lto \Ext_{\catc}^{\bdot}(\calm,\calm)\quad,\quad f\mapsto f_{\calm}\,,
\end{align*}
is a homomorphism of graded rings with image in the graded centre of $\Ext_{\catc}^{\bdot}(\calm,\calm)$. Thus, $ev_{\calm}$ endows $\Ext_{\catc}^{\bdot}(\calm,\calm)$ with the structure of a graded $ \fZ_{gr}^{\bdot}(\catc)$--algebra. 

Moreover, for each pair of objects
$\calm,\caln$ in $\catc$, the group $\Ext^{\bdot}_{\catc}(\calm, \caln)$ is a graded bimodule, as a graded right module over $\Ext^{\bdot}_{\catc}(\calm, \calm)$ and a graded left module over $\Ext^{\bdot}_{\catc}(\caln, \caln)$. It becomes therefore a graded bimodule over $ \fZ_{gr}^{\bdot}(\catc)$ via $ev_{\calm}$ on the right and $ev_{\caln}$ on the left. These structures anticommute, in that
\begin{align*}
ev_{\caln}(f)\cdot\vp &= T^{m}(f_{\caln})\circ\vp = (-1)^{mn}f_{T^{m}\caln}\circ\vp
= (-1)^{nm}\vp\circ f_{\calm}\\ &= (-1)^{nm}\vp\cdot ev_{\calm}(f)
\end{align*}
for $f\in  \fZ_{gr}^{n}(\catc), (\vp:\calm\to T^{m}\caln)\in \Ext^{m}_{\catc}(\calm, \caln)$.
Thus, $\Ext^{\bdot}_{\catc}(\calm, \caln)$ is naturally a {\em graded symmetric bimodule\/} over
$ \fZ_{gr}^{\bdot}(\catc)$ and the Yoneda pairing on the $\Ext$--groups is {\em bilinear\/} over the graded centre of $\catc$. 
\end{sit}

\subsection{The characteristic homomorphism}
We first recall the setup of Fourier-Mukai transformations; see, for example, \cite{Cal1}.
\bdfn
Let $X\to Y, X'\to Y$ be morphisms of complex spaces, with $p:X\times_{Y} X' \to X, p':X\times_{Y} X'\to X'$ the projections from the fiber product. Let
$\calf\in D(X\times_{Y} X')$ be a complex whose cohomology is supported on $Z\subseteq X\times_{Y} X'$ with $p'|Z:Z\to X'$ finite. These data define the {\em Fourier-Mukai transformation\/}, with kernel $\calf$,
\begin{align*}
\Phi_{\calf} = p'_{*}\left(\bbbl p^{*}(?)\dotimes_{\calo_{X\times_{Y} X'}}\calf\right): D(X)\lto D(X')\,,
\end{align*}
an exact functor between the indicated triangulated categories.

A  morphism $(a:\calf\to T^{m}\calg)\in \Ext^{m}_{X\times_{Y} X'}(\calf,\calg)$, 
with the supports of $\calf,\calg$ finite over $X'$, defines a morphism between the corresponding Fourier-Mukai transformations,
\begin{align*}
\Phi_{a} := p'_{*}\left(p^{*}(?)\dotimes_{\calo_{X\times_{Y} X'}}a\right):\Phi_{\calf}\lto \Phi_{T^{m}\calg}\cong T^{m}\Phi_{\calg}\,.
\end{align*}
Remark that there is no need here to derive the functor $p'_{*}$, as it is already exact on the subcategory of those complexes whose cohomology has finite support over $X'$.
\edfn

For $\calf = \calox$, the structure sheaf of the diagonal in $X\times_{Y} X$, the associated 
Fourier-Mukai transformation is the identity on $D(X)$, so $\Phi_{\calox} = \id_{D(X)}$. 
A morphism $(g:\calox\to T^{m}\calox)\in \Ext^{m}_{X\times_{Y} X}(\calox,\calox)$ yields 
consequently a morphism $\Phi_{g}: \id_{D(X)}\to T^{m}$ of endofuntors on $D(X)$, and the commutation rules for $\dotimes$ and $T$ imply that indeed $\Phi_{g}\in \fZ^{m}_{gr}(D(X))$. 
This establishes the following result. 

\begin{prop}
\label{h.20}
For any morphism $f:X\to Y$ of complex spaces, there exists a natural homomorphism of 
graded algebras
$$
\eta_{X/Y}\cong \Phi_{?}: \Ext^{\bdot}_{X\times_{Y} X}(\calox,\calox)
\lto \fZ_{gr}^{\bdot}(D(X))\,.
$$
\qed
\end{prop}

If $g:Y\to Z$ is a morphism of complex spaces, then $\eta_{X/Y}$ factors as
\begin{align*}
\eta_{X/Y}:\Ext^{\bdot}_{X\times_{Y} X}(\calox,\calox)\xto{\rho_{Y/Z}}\Ext^{\bdot}_{X\times_{Z} X}(\calox,\calox)\xto{\eta_{X/Z}}\fZ_{gr}^{\bdot}(D(X))
\end{align*}
with $\rho_{Y/Z}$ induced from the natural morphism $X\times_{Y} X\to X\times_{Z} X$ defined by $f$ and $g$. It would thus be enough to define $\eta_{X/Z}$ for $Z$ a point and then to compose with $\rho_{Y/Z}$ for the general case.

\brem
If $g:Y\to Z$ is a morphism of complex spaces such that $gf:X\to Z$ is {\em flat}, for example, taking $Z$ to be a point, then there is the commutative diagram of homomorphisms of graded algebras
\begin{align*}
\bdi
\Ext^{\bdot}_{X\times_{Y} X}(\calox,\calox)&\rTo^{\rho_{Y/Z}}&\Ext^{\bdot}_{X\times_{Z} X}(\calox,\calox)
&\rTo^{\eta_{X/Z}}&\fZ_{gr}^{\bdot}(D(X))\\
\dTo&&\dTo_{\cong}\\
\HH^{\bdot}_{X/Y}(\calox)&\rTo&\HH^{\bdot}_{X/Z}(\calox)
\edi
\end{align*}
with the vertical comparison maps coming from \ref{h.17.5} and the algebra homomorphism at the bottom from \ref{h.23}. According to \ref{h.19}, the vertical morphism on the right is an isomorphism, whence we obtain, in a roundabout way, homomorphisms of graded {\em commutative\/} algebras
\begin{align*}
\HH^{\bdot}_{X/Y}(\calox)\lto \fZ_{gr}^{\bdot}(D(X))
\end{align*}
\erem

We now show directly that, without any flatness assumptions, the morphism $\eta_{X/Y}$ above factors through the comparison homomorphism $\Ext^{\bdot}_{X\times_{Y} X}(\calox,\calox)\to \HH^{\bdot}_{X/Y}(\calox)$, thus, that Hochschild cohomology is closer to the graded centre of the derived category
than the self-extensions of the diagonal in $X\times_{Y}X$. The key, as always, is to replace the fibre product $X\times_{Y}X$ by its derived version by means of extended resolvents. This argument applies as well to arbitrary morphisms of schemes, thus, covers, for example, schemes over the integers, when there is no recourse to a flat situation.

\begin{theorem}\label{h.20.1}
For a morphism $f:X\to Y$ of complex spaces, there exists a natural homomorphism of graded commutative algebras
$$
\chi_{X/Y}\cong \Phi_{?}:\HH^{\bdot}_{X/Y}(\calox)
\lto \fZ_{gr}^{\bdot}(D(X))
$$
that factors the map $\eta_{X/Y}$ in \ref{h.20} through the comparison map from  \ref{h.17.5}.
\end{theorem}

The {\em Proof of Theorem \ref{h.20.1}\/}, including an analysis of the structure of this characteristic homomorphism,  will occupy the remainder of this section.

We work with an extended resolvent  $\fX=(X_*,W_*,\calrs,\calbs)$ of $X/Y$ as considered in \ref{h.5a}. Recall from \ref{h.5} that $\calbs$ represents a 
locally free (or even free) DG $\calss$-algebra resolution of $\calrs$ over
$\calss=\calrs\tilde\otimes_{\caloy}\calrs$. The algebra homomorphisms
\begin{align}
j_{1}&:\calrs\xto{\cong}\calrs\otimes 1\subseteq \calss\quad\text{and}
\label{eq:j1}\\
j_{2}&:\calrs\xto{\cong}1\otimes \calrs \subseteq \calss
\label{eq:j2}
\end{align}
define, by restriction of scalars, two flat $\calrs$-algebra 
structures on $\calbs$, and the corresponding structure maps, again denoted $j_{1,2}:\calrs\to \calbs$, are quasiisomorphisms of algebras. It follows that
$$
\bar\calbs:=
\calb_*\otimes_{1\otimes \calrs}\caloxs
$$
carries two algebra structures as well: an $\calrs=\calrs\otimes
1$--structure from the action on the first factor, and the, again flat, 
$\caloxs$--structure from the second factor.
Moreover, the given (quasiiso-)morphism of algebras $\nu:\calbs\to\calrs\to\caloxs$ induces a morphism
$$
\bar\nu:=\nu\otimes_{1\otimes\calrs} 1:\bar\calbs\to\caloxs\otimes_{\calrs}\caloxs
\cong \caloxs\,.
$$
The central tool is now the following lemma.

\begin{lem}\label{h.21}
With notation as just introduced, the following hold.
\begin{enumerate}
\item For any DG $\calrs$--module $\calms$, the natural morphism
$\calms \to \calms\otimes_{\calrs\otimes 1} \calbs $
is a quasiisomorphism; in particular, the functor\footnote{Quillen, 
in \cite{Qui4}, denotes this functor $\calms\mapsto \calms\otimes^{!}_{\calrs}$ 
in the affine case.}  
$\calms\mapsto\calms\otimes_{\calrs\otimes 1}\calbs$, where the target is considered a $\calrs$--module through the action of $1\otimes\calrs$ on the right, 
is exact on the category of DG $\calrs$--modules.

\item
For any DG $\caloxs$--module $\calms$, the natural morphism
$$
\calms\otimes_{\calrs\otimes 1} \bar\calbs
\xto{\id\otimes \bar\nu} \calms\otimes_{\calrs}\caloxs\cong \calms
$$
is a quasiisomorphism in each of the following two cases: 
\bnum[(i)]
\item $\calox$
is flat over $Y$, or 
\item $\calms$ is K--flat over
$Y$ in the sense of \cite{Sp}, that is, $\calms\tilde\otimes_{f^{-1}\caloy}f^{-1}(?)$ is exact on $D(Y)$.
\enum
Furthermore, the (right) DG $\caloxs$--module $\calms\otimes_{\calrs\otimes 1} \bar\calbs$ represents in each of these cases
$\calms\dotimes_{\calrs\otimes 1} \bar\calbs$ as well as
$\calms\dotimes_{\calrs\otimes 1} \calbs\dotimes_{1\otimes\calrs}\caloxs$.

\item If $\calox$ is flat over $Y$ then the functor $\calms\mapsto
\calms\otimes_{\calrs\otimes 1} \bar\calbs$ is exact on DG
$\caloxs$--modules.

In general, the functor $?\dotimes_{\calrs\otimes 1} \bar\calbs$ defines an exact auto-equivalence on $D(\caloxs)$ that is isomorphic to the identity functor via $\id\otimes\bar\nu$.
\end{enumerate}
\end{lem}

\begin{proof}
It is sufficient to show the corresponding
statements for every stalk. Thus, for a simplex $\alpha$ and
$x\in X_\alpha$ consider the stalks
$$
M:=\calm_{\alpha,x}, \quad R:=\calr_{\alpha,x}, \quad
S:=\cals_{\alpha,x},\quad
A:=\calo_{X_{\alpha,x}}, \quad  B:=\calb_{\alpha,x}\quad
\bar B:=\bar\calb_{\alpha,x}.
$$
To deduce (1) consider the DG algebra $R':=R\otimes_{R^0}B^0
\cong (R\otimes 1)\otimes_{R^0\otimes 1}B^0$.
As $R^0\otimes 1\to B^0$ is flat, the functor
$M\mapsto M\otimes_RR'\cong M\otimes_{R^0}B^0$ is exact. Moreover,
since $B$ is a free DG algebra over $R'$, by
\cite[2.17 (1)]{BFl} the functor $-\otimes_{R'}B$ is exact, whence
$M\mapsto M\otimes_RB\cong M\otimes_{R}R'\otimes_{R'}B$ is exact as
well. This proves exactness of the functor 
$\calms\mapsto\calms\otimes_{\calrs\otimes 1}\calbs$. 

It remains to show for (1) that the natural map
$M\to M\otimes_R B=M\otimes_{R\otimes 1}B$ is a quasiisomorphism.
As tensor products are
compatible with direct limits, we may suppose that $M$ is a
finitely generated $R$--module. As the functor
$-\otimes_{R\otimes 1}B$ is exact, after replacing $M$ by a
quasiisomorphic complex of free $R$--modules, we may suppose
that $M$ is free over $R$. By \cite[2.13 (2)]{BFl} the submodules
$M^{(i)}$ generated by all elements of degree $>i$ form a filtration
of $M$ by free DG submodules, and the quotients $M^{(i-1)}/M^{(i)}$
are as well free. Using a simple spectral sequence argument, we are
thus reduced to the case that $M$ is freely generated in one degree
so that, up to a shift, $M$ is isomorphic  as a DG module to a direct
sum of copies of $R$. In this case, the map in question is a direct sum of copies of the quasiisomorphism $j_{1}:R\cong R\otimes 1\to B$, and the assertion follows.

To prove the first claim in (2), it suffices, in view of (1), to show that
$M\otimes_RB\to M\otimes_R\bar B$ is a quasiisomorphism, as the composition
of the sequence of morphisms 
\begin{align*}
M\to M\otimes_{R\otimes 1}B\to M\otimes_{R\otimes 1}\bar B 
\xto{\id_{M}\otimes\bar\nu} 
M\otimes _{R} A\otimes_{R}A \cong M
\end{align*}
is the identity on $M$, and the leftmost morphism is a quasiisomorphism by part (1).
We will show this substitute claim more generally, replacing $B$ by any complex of free
$S$--modules, say, $F$ that is bounded above;
we have to replace then, of course, $\bar B$ by
$\bar F:=F\otimes_{1\otimes R} A$, and the assertion becomes that in this situation
\begin{align}
\label{claim}
\tag{$*$}
M\otimes_{R}F \lto M\otimes_{R}{\bar F}
\end{align}
is a quasiisomorphism.

Now assume first that (ii) holds. As before, using again \cite[2.13 (2)]{BFl}, the
$R$--submodules $F^{(i)}$  of $F$ that are generated by all elements of degree $>i$ 
form a filtration of $F$ by free DG submodules with free successive quotients
$F^{(i-1)}/F^{(i)}$. Again, with a simple spectral
sequence argument we reduce to the case that up to a
shift $F$ is isomorphic to a direct sum of
copies of $S$ as a DG module over $S$. Clearly we may suppose that
$F$ is of rank 1, that is, $F\cong S$.
In this case, $M\otimes_R F \cong M\otimes_R S$ is isomorphic to
the analytic tensor product $M\tilde\otimes_\Lambda R$, where
$\Lambda:=\calo_{Y,f(x)}$. As $M$ is K--flat over $\Lambda$, the functor
$M\tilde\otimes_\Lambda-$ is exact, whence $M\tilde\otimes_\Lambda R$ is
quasiisomorphic to 
$M\tilde\otimes_\Lambda A\cong M\otimes_{R\otimes 1}S\otimes_{1\otimes R}A$, 
as required.

Next assume that (i) holds in (2). $M$ can be written as a direct
limit of subcomplexes that  are bounded above.  As tensor products are
compatible with direct limits, we may suppose that
$M$ itself is bounded above. Now we can proceed as in the previous case
and reduce the assertion to the case that
$F=S$ so that as before $M\otimes_RS\cong M\tilde\otimes_\Lambda R$.
As $A$ is flat over $\Lambda$, the mapping cone over $R\to A$ is an
exact complex of $\Lambda$--flat finite $R^0$--modules. Thus it remains
exact after applying $M^p\tilde \otimes_\Lambda-$ for any $p$,
and so $M^p\tilde\otimes_\Lambda R$ is quasiisomorphic to
$M^p\tilde\otimes_\Lambda A$. Taking total complexes proves the assertion that the morphism in (\ref{claim}) is a quasiisomorphism also in this case. 

Concerning the final claim in (2), note that 
$\calbs\otimes_{1\otimes \calrs}\caloxs$ always represents 
$\calbs\dotimes_{1\otimes \calrs}\caloxs$, as $\calbs$ is flat over $R$ via $j_{2}:R\cong 1\otimes R\to B$. 
In case (i), consider the following commutative diagram
\begin{align*}
\bdi
&&B&\rTo&\bar B=B\otimes_{1\otimes R}A\\
&&\uTo_{free}&&\uTo_{free}\\
R\otimes 1&\rTo^{j_{1}}&S&\rTo&S\otimes_{1\otimes R}A\cong R\tilde\otimes_{\Lambda}A\\
&&\uTo_{j_{2}}&&\uTo\\
&&1\otimes R&\rTo&A
\edi
\end{align*}
in which the two squares, as well as the rectangle they form, represent tensor products of algebras. As $S\to B$ is flat (even free), and $\bar B \cong B\otimes_{S}(S\otimes_{1\otimes R}A)$, we obtain that $\bar B$ is flat (even free) over $S\otimes_{1\otimes R}A$. On the other hand, $A$ is flat over $\Lambda$ by hypothesis (i), and so $R\cong R\otimes 1\to R\tilde\otimes_{\Lambda}A$ is a flat homomorphism of algebras. In summary, 
this exhibits $j_{1}:R\cong R\otimes 1\to \bar B$ as the composition of the two flat maps
$R\cong R\otimes 1 \to R\tilde\otimes_{\Lambda}A$ and 
$R\tilde\otimes_{\Lambda}A\to \bar B$. The second factor being flat and bounded above, $M\otimes_{R\otimes 1}\bar B$ represents
$M\dotimes_{R\otimes 1}\bar B$. In case (ii), the first factor $M$ is K--flat by assumption, whence $M\otimes_{R\otimes 1}\bar B$ again represents  $M\dotimes_{R\otimes 1}\bar B$.

Finally, (3) is an immediate consequence of (2) and the fact that any DG module admits a K--flat resolution.
\end{proof}

\begin{rem}\label{h.22}
Note that the preceding result is already non-trivial and interesting
in the affine case. For instance, let $A=\Lambda/I$ be a quotient of a
regular local ring $\Lambda$ containing $\bbbq$, and let $R\to A$ be a resolution of $A$
by a free DG $\Lambda$--algebra; of course we may assume that $R^0=\Lambda$.
Choose an algebra resolution $B$ of $S:=R\otimes_\Lambda R\to R$ 
so that $B\to R$ is a quasiisomorphism and $B$ is free as $S$-algebra. 

Given an $A$--module $M$ we can construct a free resolution of $M$ over $A$ as follows:

Choose a free resolution $F$ of $M$ over $\Lambda$\/ that admits a DG module structure over $R$.
By \ref{h.21} (1) above, the complex
$F\otimes_{R\otimes 1} B\otimes_{1\otimes R}A$
constitutes then an $A$--free resolution of $M$.

For instance, if $A=\Lambda/f\Lambda$ is a hypersurface, then we can take $R$
to be the Koszul complex over $f$,
i.e.\ the free $\Lambda$-algebra $\Lambda[\varepsilon]$ with a generator
$\varepsilon$ of degree -1 and differential
$\partial(\varepsilon) = f$. Thus
$S=\Lambda[\varepsilon_1,\varepsilon_2]$ with
$\varepsilon_1:=1\otimes \varepsilon$ and $\varepsilon_2:=
\varepsilon\otimes 1$. The free $S$-algebra $B:=S[\eta]$
with a generator
$\eta$ of degree $-2$ and differential 
$\partial(\eta)=\varepsilon_1 -\varepsilon_2$ constitutes an
algebra resolution\footnote{This is no longer true in positive characteristic, say $p>0$, as then
$\partial(\eta^{pn})=0$. One either has to resolve further, or, more economically, use a divided power algebra.
This is precisely the place, where the (DG algebras arising from) simplicial algebras become advantageous.}
 of $R$ via the map $\eta\mapsto 0$ and
$\varepsilon_i\mapsto \varepsilon$. If now
$M$ is a maximal Cohen--Macaulay module over $A$ then its minimal
resolution
$$
0\to F_1\to F_0\to M\to 0
$$
over $\Lambda$ is of length 1. Since $M$ is annihilated by $f$, we have
$fF_0\subseteq F_1$, whence multiplication by $f$ yields a map
$\varepsilon:F_0\to F_1$. Via this map the complex
$F=(F_1\to F_0)$ has a natural DG $R$--module structure. The
resolution above is in this case just the periodic resolution
constructed by Eisenbud and Shamash
\cite{Eis, Sh}; the periodicity is given by multiplication by $\eta$.
Applying this in a similar way to complete intersections
$A=\Lambda /(f_1,\ldots, f_r)$ leads to \cite[2.4]{ABu}. 
\end{rem}

\begin{sit}
We now turn to the explicit description of the algebra homomorphism $\chi$.
To this end, we exhibit the right action of $\HH^{\bdot}_X(\calox)$ via
$ev_{\calm}\circ\chi_{X}$ on $\Ext^{\bdot}_X(\calm,\caln)$. We know from
 \cite[2.28]{BFl}, see \ref{h.4}(\ref{C.3}), and \ref{h.21} (2.(i)) that
\begin{align*}
\Ext^{\bdot}_X(\calm,\caln)&\cong \Ext^{\bdot}_{X_*}(\calms,\calns)\\
&\cong
\Ext^{\bdot}_{X_*}(\calms\otimes_{\calrs\otimes 1} \bar\calbs ,
\calns\otimes_{\calrs\otimes 1} \bar\calbs)\,,
\end{align*}
{where we consider $\calms\otimes_{\calrs\otimes 1} \bar\calbs$ and
$\calns\otimes_{\calrs\otimes 1} \bar\calbs$ as $\caloxs$-modules
from the right. We also have from \ref{h.17}(1) that}
\begin{align*}
\HH_X^{\bdot}(\calox)&\cong \Ext^{\bdot}_\calss(\calbs, \calbs)\,.
\end{align*}
The required module structure is now induced by the pairing
$$
\Ext^{\bdot}_{X_*}(\calms,\calns)\times \Ext^{\bdot}_\calss(\calbs, \calbs)
\lto \Ext^{\bdot}_{X_*}(\calms\dotimes_{\calrs\otimes 1} \bar\calbs ,
\calns\dotimes_{\calrs\otimes 1} \bar\calbs )
$$
that sends $(f,g)$ to  $f\dotimes_{\calrs\otimes 1} g\dotimes_{1\otimes \calrs}\caloxs$.
We thus obtain $\chi(g)_{\calm}$, for a homogeneous element $g\in\Ext^{|g|}_{\calss}(\calbs,\calbs)$, interpreted as a morphism $g:\calbs\to T^{|g|}\calbs$ in $D(\cals)$, by applying the \Cech functor to the morphism on top of the commutative diagram
$$
\bdi
\calms&\rTo^{\chi(g)_{\calms}}&T^{|g|}\calms\\
\uTo^{\simeq}&&\uTo_{\simeq}\\
\calms\dotimes_{\calrs\otimes 1}\bar\calbs&\rTo^{\id_{\calms}\otimes \bar g }&\calms\dotimes_{\calrs\otimes 1}T^{|g|}\bar\calbs
\edi
$$
where $\bar g =g\otimes \id_{\caloxs}$ and the vertical arrow on the left represents the quasiisomorphism from 
\ref{h.21}(2), while the vertical arrow on the right is the corresponding 
quasiisomorphism for $T^{|g|}\calm$, composed with the (quasi)isomorphism 
$$
\calms\dotimes_{\calrs\otimes 1}T^{|g|}\bar \calbs
\xto{\simeq}(T^{|g|}\calms)\dotimes_{\calrs\otimes 1}\bar\calbs
\,.
$$
One sees now immediately that for every $g\in \Ext^{m}_\calss(\calbs, \calbs)\cong 
\HH_X^{m}(\calox)$, the family 
$$
\chi(g)_{\calm} =ev_{\calm}\chi(g)
:= (?\otimes {\bar g}):\calm\to T^{m}\calm \quad \text{for $\calm \in D(X),$}
$$ 
defines an element of $\fZ_{gr}^{m}(D(X))$, and that $\chi$ is a homomorphism of graded algebras.
\end{sit}

This completes the {\em Proof of Theorem \ref{h.20.1}\/}. \qed

\begin{rems}
(1) This characteristic homomorphism to the graded centre of the derived category has made its appearence in several forms for the affine situation of algebra homomorphisms, for example in
\cite[Sect.3]{ABu}, \cite{ASu}, \cite{SSo} or \cite{ET}.

(2) For compact manifolds, C\u{a}ld\u{a}raru, in \cite{Cal1, Cal2}, investigates higher 
order structure on the Hochschild cohomology of compact complex manifolds and notes 
the existence of the characteristic homomorphism in \cite[4.10]{Cal1}. He also points 
out that for the example of  an elliptic curve $X$, the homomorphism $\chi_{X}$ is not 
injective: $HH^{2}_{X}(\calox)\cong H^{1}(X,\calox)$ is a onedimensional vector space, 
but $\fZ^{2}_{gr}(D(X)) = 0$ as the category of quasicoherent sheaves on $X$ is 
hereditary, that is, of global dimension one. 
\end{rems}

\section{Functoriality of the Hochschild complex} In this section
we will show that the Hochschild complex is a well defined object of
the derived category and that it behaves functorially  with respect to
morphisms of complex spaces. The idea of the proof is the same as the one 
for the cotangent complex given in \cite{Fle1}; however we will
formalize the treatment as it seems useful as well in other situations, such as the 
one considered here. Moreover we give the full
details of the (non trivial) proof that was left to the
reader in ({\em loc.cit.\/}).

\subsection{Categories of models}

Assume that we are given a commutative diagram of categories\footnote{We will henceforth only consider {\em small\/} categories.} and
functors as indicated by the solid arrows:
\bdi[h=7mm]
\catm & \rTo^{G}  & \catf \\
\dTo<{F} & \ruDotsto^{\bar G}& \dTo_{H} \\
\catc &\rEq & \catc \,.
\edi
In this subsection we will give a simple criterion as to when the
given functor $G:\catm\to \catf$ can be factored through $F$ by a functor
$\bar G$, represented by the dotted arrow.

To explain what we have in mind consider the following example.

\begin{exam}\label{m.1}
Let $\catc = \catmor$ be the category of holomorphic mappings of complex
spaces, where we write an object $X\to Y$ of $\catmor$ in brief as $X/Y$. Recall that the
morphisms $f=(f_X,f_Y):X'/Y'\to X/Y$ in $\catmor$ are given by
commutative diagrams
\be\label{morXY}
\bdi[h=7mm]
X'& \rTo^{f_X}  & X\\
\dTo&&\dTo\\
Y'&\rTo^{f_Y} & Y
\edi
\ee
Furthermore, let $H:\catf\to\catmor$ be the fibration in derived
categories of complexes of modules that are bounded above so that 
the fibre of $\catf\to \catmor$ over an object $X/Y$ is just the derived 
category  $D^-(X)$ and the
fibre functor over a morphism $f=(f_X,f_Y):X'/Y'\to X/Y$ is given by 
$Lf_{X}^{*}:D(X)\to D(X')$.

For the category of models $\catm$, choose the
category of extended resolvents of morphisms of complex spaces; we will define the morphisms in that  category in \ref{m.6}. Let $F:\catm\to \catc=\catmor$ be the functor that assigns to the extended resolvent the map it resolves. Finally, for every
extended resolvent $\fX=(X_*,W_*,\calrs,\calbs)$ of $X/Y$ we have defined in \ref{h.5} the
Hochschild complex with respect to this extended resolvent 
$$
G(\fX)=\bbbh_{\fX}:=\calcb((\calbs\otimes_\calss\caloxs)|X_*)\,,
$$
and we need to know that this complex does
not depend on the choice of the extended resolvent, equivalently, that there is a
functor ${\tilde G}:\catmor\to \catf$ that assigns to every morphism $X\to Y$ just one 
Hochschild complex $\bbbh_{X/Y}$.
\end{exam}

In the following, the objects of the category $\catm$ over a given
object $X$ of $\catc$ will be called {\em models} of $X$. With this
notation we introduce the following definition.

\begin{defn}\label{m.2}
The category $\catm$; or, more precisely, the functor $F:\catm\to\catc$; is called a {\em 
category of models}\/ for $\catc$ if the following conditions are satisfied.

(1) Every object $X$ in $\catc$ admits a model $\fX\in\catm$, that
is, $F$ is surjective on objects.

(2) Let $f:X\to Y$ be a morphism in $\catc$. Given
models $\fX$ of $X$ and $\fY$ of $Y$, there is a model
$\fX'$ of $X$ together with a morphism $\fX'\to \fX$ over the
identity of $X$ and a morphism $\fX'\to \fY$ over $f$.

(3) Let $f:X\to Y$ be a morphism in $\catc$. If we are given two
morphisms $\ff_1,
\ff_2:\fX\to\fY$ over $f$ then we can find a model $\fX'$ of
$X$ together with commutative diagrams of models for $i=1,2$: 
$$
\bdi[h=7mm]
&&\fX\\
& \ldTo^{\id} && \rdTo>{\ff_i}\\
\fX && \dTo>{\fg_i} && \fY\\
& \luTo<{\pi} && \ruTo>{\ff'}\\
&&\fX'
\edi
\qquad\mbox{over}\qquad
\bdi[h=7mm]
&& X\\
& \ldTo^{\id} && \rdTo>{f}\\
X && \dTo>{\id} && Y\\
& \luTo<{\id} && \ruTo>{f}\\
&&X&&.
\edi
$$
\end{defn}

\begin{sit}
Consider the set of morphisms $\Sigma$ in $\catm$ that lie over an
isomorphism of $\catc$, so that a morphism $\ff$ of $\catm$ belongs to
$\Sigma$ if and only if the underlying morphism $f = F(\ff)$ in $\catc$ is an
isomorphism. With respect to this set, one can form the {\em localized category\/}
 $\catm[\Sigma^{-1}]$ as in \cite[1.1]{GZi}. Recall that the objects of $\catm[\Sigma^{-1}]$ are
just the objects of $\catm$, whereas a morphism $\fX\to\fY$ in
$\catm[\Sigma^{-1}]$ consists of an equivalence class of a zigzag
\begin{equation}\label{zigzag}
\bdi[s=7mm]
  && \fY_1 &&&& \fY_2 &&\ldots &&\fY_n\\
&\ldTo<{\beta_1} &&\rdTo>{\alpha_1} &&
\ldTo<{\beta_2} &&  \ldots &&\ldTo<{\beta_n} &&\rdTo>{\alpha_n}\\
\fX=\fX_0 &&&& \fX_1 &&\ldots && \fX_{n-1} &&&&\fX_n=\fY\,,
\edi
\end{equation}
where $F(\beta_i)$ is an isomorphism for $i=1,\ldots, n$.
The construction comes with the obvious localization functor $p:\catm\to \catm[\Sigma^{-1}]$ that is the identity on objects and considers each morphism $\ff:\fX\to \fY$ from $\catm$ as the zigzag 
$
\bdi
\fX&\lTo^{\id}&\fX&\rTo^{\ff}&\fY
\edi
$
in $\catm[\Sigma^{-1}]$.

By the universal property of such localized categories; see \cite[1.1 (ii)]{GZi}; there is a canonical factorization of $F$ through the localization, thus, a commutative diagram
\bdi[h=7mm]
\catm &&\rTo^p && \catm[\Sigma^{-1}]\\
& \rdTo<F && \ldTo>{\bar F}\\
  && \catc
\edi
where $\bar F:\catm[\Sigma^{-1}]\to \catc$  is the induced functor.
\end{sit}

The main insight here is now contained in the following result.

\begin{theorem}\label{m.3}
If $\catm$ is a category of models for $\catc$ then the functor
$\bar F$ is an equivalence of categories.
\end{theorem}

For the proof we need the following lemma.

\begin{lem}\label{m.4}
(1) Assume given a model $\fX$ of $X$, a model $\fY$ of $Y$, and a morphism 
$f:X\to Y$ in $\catc$. If $\ff_1,\, \ff_2:\fX\to\fY$ are
two morphisms over $f$ then the localized morphisms $p(\ff_1),\,
p(\ff_2): p(\fX)\to p(\fY)$ in $\catm[\Sigma^{-1}]$ are equal.

(2) Let $\bar h:p(\fX)\to p(\fY)$ be a morphism in
$\catm[\Sigma^{-1}]$ over $f:X\to Y$ in $\catc$. There is then a model
$\fX'$ of $X$ and morphisms $\fg:\fX'\to\fX$ over $\id_X$ and
$\ff:\fX'\to \fY$ over $f$ such that $\bar h\circ p(\fg)=p(\ff)$.
\end{lem}

\begin{proof}
(1) If we choose a model $\fX'$ of $X$ as in \ref{m.2} (3) then
$p(\pi)$ is an isomorphism in $\catm[\Sigma^{-1}]$ whence
$p(\fg_1)=p(\pi)^{-1}=p(\fg_2)$ and then also  $p(\ff_1)=p(\ff_2)$.

For (2), it suffices to show by \ref{m.2}(2) that every morphism in
$\catm[\Sigma^{-1}]$ is represented by a zigzag (\ref{zigzag}) of
length $n=1$. Assume that the morphism $\bar h$ is
represented by a zigzag as in (\ref{zigzag}), of minimal length $n\ge
1$. If $n>1$ then we can find a model $\fY_{n-1}'$ of
$Y_{n-1}:=F(\fY_{n-1})$ together with two morphisms
$\fg_{n-1}:\fY_{n-1}'\to \fY_{n-1}$ over $\id_{Y_{n-1}}$ and
$\fg_n:\fY_{n-1}'\to \fY_n$ over $F(\beta_n)^{-1}\circ
F(\alpha_{n-1})$, that is, the diagram 
\bdi[s=7mm]
F(\fY_{n-1}')&\rTo^{F(\fg_n)}&F(\fY_n)\\
\dTo<{F(\fg_{n-1})=\id} &&  \dTo>{F(\beta_n)}\\
F(\fY_{n-1}) & \rTo>{F(\alpha_{n-1})}& F(\fX_{n-1})
\edi
commutes. Using part (1), the diagram
\bdi[s=7mm]
p(\fY_{n-1}')&\rTo^{p(\fg_n)}&p(\fY_n)\\
\dTo<{p(\fg_{n-1})} &&  \dTo>{p(\beta_n)}\\
p(\fY_{n-1}) & \rTo>{p(\alpha_{n-1})}& p(\fX_{n-1})
\edi
already commutes in $\catm[\Sigma^{-1}]$, whence $\bar h$ is represented
by the shorter zigzag
\bdi[s=7mm]
  && \fY_1 &&&& \fY_2 &&\ldots &&\fY_{n-2}&&&&\fY_{n-1}'\\
&\ldTo<{\beta_1} &&\rdTo>{\alpha_1} &&
\ldTo<{\beta_2} &&  \ldots&&\ldots &&\rdTo
&&\ldTo<{\beta_{n-1}'} &&\rdTo>{\alpha_{n-1}'}\\
\fX=\fX_0 &&&& \fX_1 &&\ldots &&\ldots &&&& \fX_{n-2} &&&&\fX_n=\fY\,,
\edi
where $\beta_{n-1}':=\beta_{n-1}\circ\fg_{n-1}$ and
$\alpha_{n-1}':=\alpha_n\circ\fg_n$.
\end{proof}

\begin{proof}[Proof of theorem \ref{m.3}]
As $F$ is surjective on objects, so is $\bar F$. In order to show that $\bar F$
is fully faithful, let $\fX$, $\fY$ be objects of $\catm$ over $X$, resp.\ $Y$, 
and consider the maps
$$
\Mor_\catc (\fX, \fY)\xto{p}\Mor_{\catm[\Sigma^{-1}]}(p(\fX), p(\fY))
\xto{\bar F}\Mor_\catc(X,Y)\,.\leqno (*)
$$
We need to show that the second map, labeled ${\bar F}$ is bijective.
We first show that it is surjective. In view of condition \ref{m.2} (2), for a given
morphism $f:X\to Y$, we can find a model
$\fX'$ of $X$ together with morphisms $\fg:\fX'\to \fX$ over $\id_X$ and
$\ff:\fX'\to\fY$ over $f$, whence the morphism $p(\ff)\circ
p(\fg)^{-1}\in \Mor_{\catm[\Sigma^{-1}]}(p(\fX), p(\fY))$ maps to $f$ under 
$\bar F$.

To prove the injectivity of the second map in $(*)$, let $\bar
f_1\,,\bar f_2: p(\fX)\to p(\fY)$ be morphisms over the same morphism
$f:X\to Y$. By \ref{m.4} (2), $\bar f_1$ and $\bar f_2$ can be represented by
zigzags of length 1,
$$
\bdi
&& \fX_1\\
&\ldTo<{\fg_1} && \rdTo>{\ff_1}\\
\fX &&\rDotsto^{\bar f_1} && \fY
\edi
\qquad\mbox{and}\qquad
\bdi
&& \fX_2\\
&\ldTo<{\fg_2} && \rdTo>{\ff_2}\\
\fX &&\rDotsto^{\bar f_2} && \fY,
\edi
$$
where $\ff_i$, for $i=1,2$,  is a morphism over $f$ and $\fg_i$ is a 
morphism over $\id_X$. (in these diagrams we use the convention that the
solid arrows are morphisms in $\catm$, whereas the dotted ones are only morphisms in
$\catm[\Sigma^{-1}]$). By condition \ref{m.2} (2), we can find a model
$\tilde\fX$ and morphisms
$\fh_i:\tilde\fX\to \fX_i$, $i=1,2$, over $\id_X$.
By construction, we have then
$$
F(\fg_1\circ\fh_1)=F(\fg_2\circ\fh_2)\quad \mbox{and}\quad
F(\ff_1\circ\fh_1)=F(\ff_2\circ\fh_2).
$$
In view of \ref{m.4} (1), this implies
$$
p(\fg_1\circ\fh_1)=p(\fg_2\circ\fh_2)\quad \mbox{and}\quad
p(\ff_1\circ\fh_1)=p(\ff_2\circ\fh_2)
$$
in $\catm[\Sigma^{-1}]$. It follows now from
$$
\bar f_i = p(\ff_i)\circ
p(\fg_i)^{-1}=p(\ff_i\circ\fh_i)\circ p(\fg_i\circ\fh_i)^{-1}\,,
\quad i=1,2,
$$
that ${\bar f}_{1}={\bar f}_{2}$ as we had to show.
\end{proof}

Returning to the problem discussed at the beginning of this section we
get the following application.

\begin{cor}\label{m.5}
Let $F: \catm\to \catc$ be a category of models for $\catc$ and assume that
$G:\catm\to \catf$ is a functor such that for every morphism
$\ff:\fX\to\fY$ over an isomorphism $f:X\to Y$ in $\catc$, the morphism
$G(\ff):G(\fX)\to G(\fY)$ is also an isomorphism. Then there is a
unique functor $\bar G:\catc\to\catf$ such that $G= \bar G\circ F$.
\end{cor}
The proof follows
immediately from the the main result \ref{m.3} above and the universal
property of localizations; see \cite[1.1 (b)]{GZi}).\qed

\subsection{Resolvents as a category of models}

In a first step we introduce morphisms of resolvents.

\begin{defn}\label{m.6}
Assume we are given a morphism $f=(f_X,f_Y):X'/Y'\to X/Y$ in
$\catmor$  and extended resolvents 
$$
\fX^{\prime(e)}=(X'_*,W'_*,\calr_*',\calb_*') \mbox{ of }X'/Y'
\quad\mbox{ and }\quad
\fX^{(e)}=(X_*,W_*,\calrs,\calbs) \mbox{ of }X/Y\,.
$$
A morphism $\ff:\fX^{\prime(e)}\to \fX^{(e)}$ over $f$ consists
of the following data.
\begin{enumerate}
\item A map $\sigma:I'\to I$
with $f_X(X'_i)\subseteq X_{\sigma(i)}$; this induces a
map, denoted by the same symbol, $\sigma:A'\to A$ of simplicial sets, 
and a morphism $f_{X_*}:X'_*\to X_*$ of simplicial spaces.

\item A  (simplicial) morphism $ f_{W_*}:
W'_*\to W_*$ restricting to $f_{X_*}$ on $X_*$ such that the diagram
\bdi
W'_* & \rTo^{ f_{W_*}} & W_*\\
\dTo & & \dTo\\
Y' &\rTo & Y
\edi
commutes.

\item A morphism of $\calo_{W'_*}$--algebras
$$
{\ff}^{*}(\calr_*):=f_{W_*}^{*}(\calr_*)\to \calr_*'
$$
compatible with the projections to $\calo_{X'_*}$.

\item With $\cals_*'$ and $\calss$ as in \ref{h.2} and denoting the
pullback of an $\calo_{W_*\times_YW_*}$-algebra $\calds$ under
$f_{W_*}\times f_{W_*}:W'_*\times_YW'_*\to W_*\times_YW_*$ by
$\ff^{*}(\calds)$,  there is given a morphism
$\ff^{*}(\calbs)\to\calb_*'$ such that the diagram
\bdi
\ff^{*}(\calss)&\rTo &\ff^{*}(\calbs)&\rTo
& \ff^{*}(\calrs)\\
\dTo && \dTo && \dTo \\
\cals_*'&\rTo & \calb_*'&\rTo
& \calr_*'\\
\edi
commutes, where the first vertical arrow is the tensor
product of the morphism $\ff^{*}(\calrs)\to \calr_*'$ in (3) with itself.
\end{enumerate}
\end{defn}
One can compose morphisms of extended resolvents in an obvious way,
and the identity is a morphism of extended resolvents. Thus the
extended resolvents form a category over the category
$\catmor$ of morphisms of complex spaces.

In a similar way one can also form the category of resolvents;
morphisms of resolvents will be given by the data (1)-(3) in
\ref{m.6}.

\begin{prop}\label{m.7}
The category of (extended) resolvents constitutes a model
category over $\catmor$ in the sense of \ref{m.2}.
\end{prop}

For the proof we need a few preparations. Let $\fX$ be a resolvent of
$X/Y$ and $\calss$ a DG $\calo_{W_*}$--algebra that,
according to our conventions, is assumed to be concentrated in degrees
$\le 0$ and to have coherent homogeneous components.  By Proposition
2.19 in \cite{BFl}, the category of those $\calss$--modules that are bounded above 
with coherent cohomology has enough projectives. More precisely, for a
simplex $\alpha$ and an $\cals_\alpha$--module $\calp_\alpha$, one can form the
module $p_\alpha^{*}(\calp_\alpha)$ defined by
$p_\alpha^{*}(\calp_\alpha)_\beta:=0$, if $\alpha\not\subseteq \beta$,
and $p_\alpha^{*}(\calp_\alpha)_\beta:=
p_{\alpha\beta}^{*}(\calp_\alpha)
\otimes_{p_{\alpha\beta}^{*}(\cals_\alpha)}\cals_\beta$
otherwise, where
$p_{\alpha\beta}:W_\beta\to W_\alpha$ are the transition maps. This is
a simplicial module with respect to the obvious  transition maps,
and it is projective if $\calp_\alpha$ is projective. Moreover, as
is shown in \cite[2.13 (1)]{BFl}, every projective module
over $\calss$ is a direct sum of such modules. By definition, the
graded free modules over $\calss$ are those that admit a direct sum
decomposition $\calps=\bigoplus_\alpha p_\alpha^{*}(\calp_\alpha)$, where
each $\calp_\alpha$ is graded free over $\cals_\alpha$. We need the
following observation.

\begin{lem}\label{m.8}
Let $\ff: \fX'\to \fX$ be a morphism of resolvents and
assume that
$\calss$ is a DG $\calo_{W_*}$--algebra. If, with the notation of
\ref{m.6}, the map $\sigma:I'\to I$ is injective then the following
hold.

1. For every projective (graded free) $\calss$--module $\calps$, its
pull back $\ff^{*}(\calps)$ is again projective (graded free) over $\ff^{*}(\calss)$.

2. For every graded free $\calss$--algebra $\calts$, its pull back
$\ff^{*}(\calts)$ is again a graded free algebra over $\ff^{*}(\calss)$.
\end{lem}

\begin{proof}
(2) is an immediate consequence of (1). For the proof of (1), we may assume,
in view of the structure theorem for projective modules mentioned above,
that $\calps\cong p_\alpha^{*}(\calp_\alpha)$ for some
projective (resp.\ graded free) module over $\cals_\alpha$. If now
$\alpha\not\subseteq \sigma(I')$, then $\ff^{*}(\calps)=0$ and the
assertion holds trivially true. Otherwise $\alpha=\sigma(\alpha')$ for
a unique simplex $\alpha'$ of
$A'$, and then $\ff^{*}(\calps)\cong
p_{\alpha'}^{*}((\ff^{*}(\calp)_{\alpha'})$, whence this module is
projective, resp.\ graded free over $\ff^{*}(\calss)$.
\end{proof}

We remark that the lemma is in general no longer true if the map 
$\sigma$ fails to be injective.

\begin{proof}[Proof of \ref{m.7}]
We have to verify that the conditions (1)--(3) in \ref{m.2} are
satisfied. The existence of resolvents follows from  
\cite[2.34]{BFl}, and the existence of extended ones from the
discussion in \ref{h.5}. In order to deduce (2), let
$f=(f_X,f_Y):X'/Y'\to X/Y$ be a morphism in
$\catmor$ and let
$\fX^{\prime(e)}=(X'_*,W'_*,\calr'_*,\calb_*')$ and
$\fX^{(e)}=(X_*,W_*,\calrs,\calbs)$ be extended resolvents of
$X'/Y'$, resp.\ $X/Y$. The sets
$$
\tilde X'_k:=X'_j\cap f_X^{-1}(X_i)\,,\quad
k:=(i,j)\in \tilde I':=I\times I',
$$
form again a covering of $X'$ by Stein compact sets with associated
nerve, say $\tilde X'_*:=(\tilde X'_\alpha)_{\alpha\in \tilde
A'}$. The projections
$p:\tilde I'\to I'$ and $q:\tilde I'\to I$ provide maps of
simplicial schemes, again denoted $p:\tilde A'\to A'$ and
$q:\tilde A'\to A$, and  $\tilde X'_*$ can be embedded diagonally into
the simplicial space $\tilde W'_*$ with 
$$
\tilde W'_\alpha:=W'_{p(\alpha)}\times_Y W_{q(\alpha)}, \quad
\alpha\in \tilde A'\,,
$$
that is smooth over $Y'$.
The projections $\pi: \tilde W'_*\to W_*$ and $\pi': \tilde
W'_*\to W'_*$ are then morphisms of simplicial spaces that restrict to
the morphisms $f_X$, resp.\ $\id$, on $\tilde X'_\alpha$ for every
simplex $\alpha\in \tilde A'$. We can now choose a free algebra
resolution
$\tilde \calr_*'$ of the induced algebra homomorphism 
$$
\pi^{\prime *}(\calr'_*)\otimes_{\calo_{\tilde
W'_*}}\pi^{*}(\calrs)
\lto \calo_{\tilde X'_*}\,.
$$
The construction so far provides a resolvent
$\tilde\fX':=(\tilde X'_*, \tilde W'_*, \tilde \calr'_*)$ of
$X'/Y'$ together with morphisms of resolvents
$$
\ff':\tilde\fX':=(\tilde X'_*, \tilde W'_*,\tilde \calr'_*)\to
\fX'=( X'_*,  W'_*, \calr'_*)
\quad\mbox{and}\quad
\ff:\tilde\fX'\to \fX=( X_*,  W_*, \calrs)
$$
lying over $\id$ and $f$, respectively. To construct as well a
morphism of extended resolvents we note that $\ff'$ and $\ff$ induce
algebra homomorphisms
$$
\ff'^{*}(\cals'_*)\to
\tilde\cals'_*
\quad\mbox{and}\quad
\ff^{*}(\cals_*)\to
\tilde\cals'_*
$$
on $\tilde W_*'\times_{Y'}\tilde W'_*$ compatible with the
projections onto $\tilde\calr'_*$,  where $\tilde\cals'_*$,
$\cals_*'$,
$\calss$ are as explained in
\ref{h.2}. Hence, if $\tilde\calb_*'$ is a free algebra resolution
of 
$$
\cald_*:=\ff'^{*}(\calb'_*)\otimes_{\ff'^{*}(\cals'_*)}
\tilde\cals'_*
\otimes_{\ff^{*}(\cals_*)}\ff^{*}(\calb_*)
\lto \tilde\calr'_*
$$
on $\tilde W_*'$, then there are induced maps
$\ff'^{*}(\calb'_*)\to\tilde \calb_*'$ and
$\ff^{*}(\calb_*)\to\tilde \calb_*'$ as required in
\ref{m.2}(4), and so $\ff'$ and $\ff$ extend to morphisms of extended
resolvents $\tilde\fX^{\prime(e)}:=(\tilde X'_*,\tilde W'_*,\tilde
\calr'_*,\tilde \calb'_*)\to \fX^{\prime(e)}$  and
$\tilde\fX^{\prime(e)}\to\fX^{(e)}$, proving (2). 

In order to show that \ref{m.2} (3) is satisfied, let $\ff_i:
\fX^{\prime(e)}\to \fX^{(e)}$, $i=1,2$,  be two morphisms over
$f$ and let $f_{W_*,i}:W_*'\to W_*$ be the associated morphisms of
smoothings. With
$\tilde\fX^{\prime(e)}$ as before, the maps $\ff_i$ give rise to
morphisms of simplicial spaces
$$
g_{*i}:=(\id, f_{W_*,i}): W'_*\to \tilde W'_*=W'_*\times_Y W_*
$$
and to morphisms of simplicial algebras
$$
g_{*i}^{*}(\pi^{\prime *}(\calr'_*)\otimes_{\calo_{\tilde
W'_*}}\pi^{*}(\calrs))\cong
\calr'_*\otimes_{\calo_{W'_*}}\ff_i^{*}(\calrs))
\to \calr'_*.\leqno (*)
$$
The underlying map of simplicial schemes
$A'\to \tilde A'$ is the graph of a map and thus is injective.
Applying \ref{m.8}, $g_{*i}^{*}(\tilde \calr'_*)$ is a free algebra
over the left hand side of $(*)$ and so the morphism in $(*)$
extends to a morphism $g_{i*}^{*}(\tilde\calr'_*)\to \calr'_*$. Thus we
have constructed morphisms $\fg_i:\fX'\to\tilde \fX'$, $i=1,2$, of
resolvents as required in
\ref{m.2} (3). 

To construct as well morphism of extended
resolvents, again denoted
$\fg_i:\fX^{\prime(e)}\to\tilde\fX^{\prime(e)}$, for $i=1,2$, we note
that the data so far constructed provide for $i=1,2$ a diagram in
solid arrows 
\bdi
\calb'_*\otimes_{\cals'_*}
\fg_i^{*}(\tilde\cals'_*)
\otimes_{\ff_i^{*}(\cals_*)}\ff_i^{*}(\calb_*)
&\cong\, &
\fg_i^{*}(\cald_*)
&\rTo & \calb'_*\\
&& \dTo &\ruDotsto & \dTo\\
&&\fg_i^{*}(\tilde \calb'_*)&\rTo & \calr'_*
\edi
Applying the same argument as before, $\fg_i^{*}(\tilde \calb'_*)$ is a
free algebra over  $\fg_i^{*}(\cald_*)$, whence there is a morphism
of DG algebras $\fg_i^{*}(\tilde \calb'_*) \to\calb_*'$ as indicated
by the dotted arrow, so that the resulting diagram is still commutative. This gives the
required morphisms of extended resolvents.
\end{proof}

We are now able to deduce that Hochschild cohomology is well defined
and functorial as claimed earlier in \ref{h.main}.

\begin{theorem}\label{m.9}
(1) Assigning to an extended resolvent
$\ff:\fX^{(e)}=(X_*,W_*,\calrs,\calbs)$ the Hochschild complex
$\bbbh_{\fX^{(e)}}:=\calcb(\calbs\otimes_\calss\caloxs)\in D(X)$ defines
a functor on the category of extended resolvents. If 
$\fX^{\prime(e)}\to \fX^{(e)}$ is a morphism of extended resolvents
over an isomorphism $f=(f_X,f_Y):X'/Y'\to X/Y$ then the induced
morphism $f^{*}(\bbbh_{\fX^{(e)}})\to \bbbh_{\fX^{\prime(e)}}$ is a
quasiisomorphism. 

(2) The Hochschild complex $\bbbh_{X/Y}$ is well defined and
functorial in $X/Y$, that is, every diagram of complex spaces as displayed in
\ref{m.1} induces a well defined and functorial morphism of
algebra objects $Lf^{*}(\bbbh_{X/Y})\to \bbbh_{X'/Y'}$.
\end{theorem}

\begin{proof}
(2) is a consequence of (1) and \ref{m.5} applied to the category
$\catm$ of extended resolvents. To show (1),  let
$\ff^{*}(\calbs)\to \calb'_*$ be the morphism of DG algebras belonging
to $\ff$, see \ref{m.6}(4). This yields first a
morphism of complexes $f_{X_*}^{*}(\calbs\otimes_\calss\caloxs)\to
\calb_*'\otimes_{\cals_*'}\calo_{X'_*}$ on $X'_*$, and then, after applying the 
\Cech functor, a morphism of algebra objects
$$
f^{*}(\bbbh_\fX)=f^{*}(\calc^\sbullet(\calbs\otimes_\calss\caloxs))
\to \calc^\sbullet(f_{X_*}^{*}(\calbs\otimes_\calss\caloxs))\to
\calc^\sbullet(\calb_*'\otimes_{\cals_*'}\calo_{X_*'})=\bbbh_{\fX'}
$$
on $X'$. As $\bbbh_\fX$ is a complex of flat $\calox$-modules, the
term on the left represents $Lf^{*}(\bbbh_\fX)$. This construction is
compatible with compositions and transforms the identity into the
identity, thus proving functoriality of $\fX^{(e)}\mapsto 
\bbbh_{\fX^{(e)}}$.

The second part of (1)  can be deduced with the same reasoning as in
\ref{h.8}(2); we leave the details to the reader.
\end{proof}

\begin{cor}\label{m.10}
For every $\calox$--module $\calm$ there are natural maps
$\HH_{\bdot}^{X/Y}(\calm)\to \HH_{\bdot}^{X'/Y'}(Lf^{*}\calm)$. The map
$\HH_{\bdot}^{X/Y}(\calox)\to \HH_{\bdot}^{X'/Y'}(\calo_{X'})$ is compatible with
the algebra structures on both sides.
\end{cor}

The following fact was already used in Section 2, see \ref{h.23}.

\begin{prop}\label{m.11}
The graded algebra structure on $\HH_{X/Y}^{\bdot}(\calox)$ is
independent of the choice of the extended resolvent, and, for every
morphism
$f=(\id_X, f_Y): X/Y'\to X/Y$, the induced map
$\HH_{X/Y'}^{\bdot}(\calox)\to\HH_{X/Y}^{\bdot}(\calox)$ 
is a homomorphism of algebras.
\end{prop}

Before giving the proof we make following preparation.

\begin{sit}\label{m.12}
Let $\fX^{(e)}=(X_*,W_*,\calrs,\calbs)$ be an extended resolvent of a
morphism of complex spaces $X/Y$. Restricting the sheaf
$\calo_{W_*}$ to $X_*$ topologically gives a new smoothing
$W^\tau_*=(|X_*|,
\calo_{W_*}|X_*)$, and the restrictions $\calr^\tau_*:=\calrs|X_*$
and $\calb_*^\tau:=\calbs|X_*$ resolve
$\caloxs$ again as locally free DG algebras over $\calo_{W^\tau_*}$,
resp.\ $\cals_*^\tau=\calss|X_*$. In the following we will call
$\fX^{(e)}_\tau:=(X_*,\calo_{W^\tau_*},\calr^\tau_*,\calb_*^\tau)$ in
brief {\em the topological reduction} of $\fX$, and we will call
$\fX^{(e)}$  {\em topologically reduced} if
$\fX^{(e)}=\fX^{(e)}_\tau$. Obviously, forming the topological
reduction is compatible with morphisms of resolvents. The
natural morphism
$$
\Ext_{\calss}^{\bdot}(\caloxs, \caloxs)\lto
\Ext_{\cals^\tau_*}^{\bdot}(\caloxs, \caloxs)
$$
is compatible with the algebra structures on both sides and,
furthermore, it is an isomorphism as both sides represent the
Hochschild cohomology $\HH_{X/Y}^{\bdot}(\calox)$ in view of \ref{h.17}(1). 
\end{sit}

\begin{proof}[Proof of \ref{m.11}]
Let $\ff:\fX^{\prime(e)}=(X'_*, W'_*, \calr_*',\calb_*')\to
\fX^{(e)}=(X_*, W_*, \calrs,\calbs)$ be a morphism of resolvents over
$X/Y'\to X/Y$. By the
preceding remark, we may assume that $\fX^{\prime(e)}$ and
$\fX^{(e)}$ are topologically reduced. Let
$\sigma:A'\to A$, as in \ref{m.6},  be  the map of simplices associated
to the simplicial morphism $X'_*=(X'_\alpha)_{\alpha\in A'}\to
X_*=(X_\beta)_{\beta\in A}$. Forming the topological restrictions
$$
\calo_{\tilde W_\alpha}:=\calo_{W_{\sigma(\alpha)}}|X'_\alpha\,,
\quad
\tilde\calr_\alpha:=\calr_{\sigma(\alpha)}|X'_\alpha\,,
\quad
\tilde\cals_\alpha:=\cals_{\sigma(\alpha)}|X'_\alpha
\mbox{ and }
\tilde\calb_\alpha:=\calb_{\sigma(\alpha)}|X'_\alpha
$$
gives a new extended resolvent $\tilde\fX^{(e)}:=(X'_*,\calo_{\tilde
W_*},\tilde\calr_*,\tilde\calbs)$ of $X/Y$. Using the fact that
$\fX^{\prime(e)}$ is topologically reduced, there is a factorization
$$
\fX^{\prime(e)}\to \tilde\fX^{(e)}\to \fX^{(e)}.
$$
Restricting scalars to $\tilde\calss$ gives a morphism of
Ext--algebras
$$
\HH_{X/Y'}^{\bdot}(\calox)\cong\Ext^{\bdot}_{\cals'_*}(\calo_{X'_*},\calo_{X'_*})\lto
\Ext^{\bdot}_{\tilde\cals_*}(\calo_{X'_*},\calo_{X'_*})
\cong \HH_{X/Y}^{\bdot}(\calox)
$$
as desired.
\end{proof}

\begin{rem}\label{m.13}
1.  Let $X/Y$ be a morphism of complex spaces and $\fX$ as before a
resolvent.  With the same arguments as in \ref{m.9} it follows that
the cotangent complex
$$
\bbbl_{X/Y}:=\calc^\sbullet (\Omega^1_{\calrs/Y}\otimes_\calrs\caloxs)
$$
is a well defined object of the derived category $D(X)$ and that it
is functorial with respect to morphisms.

2. Categories of models can be used to derive (additive or non-additive)
functors. Let $\cata$ be an abelian category with enough projectives
and let
$\catm$ be the category of all pairs $(P,M)$, where $P$ is a
projective resolution of $M\in \cata$, that is, there is an exact
sequence $\cdots\to P^i\to P^{i+1}\to\cdots P^0\to M\to 0$ with each $P^{i}$
projective. Morphisms of such resolutions are defined in the usual way. 
The reader may easily verify that the functor of opposite categories
$F:\catm^\circ\to\cata^\circ$ with $(P,M)\mapsto M$ is a model
category in the sense of
\ref{m.2}. Let now $T:\catm\to \catf$ be a (not necessarily
additive!) functor into an arbitrary category $\catf$ such that for
$f:(P,M)\to (Q,N)$ the morphism $T(f)$ is an isomorphism whenever
the morphism $M\to N$ induced by $f$ is an isomorphism. Applying
\ref{m.5}, such a functor  admits a factorization $\bar T:\cata\to
\catf$. 

For instance, if $T(P,M):=H^i(\Hom(P,N))$ with a fixed object
$N\in \cata$ then $\bar T(M)$ is just the group $\Ext^i(M,N)$. 
Similarly, if there is an internal tensor product $\otimes$ on $\cata$ with 
the usual
properties and $T(P,M):=H^i(P\otimes N$ then we recover the
functors $\Tor_i(M,N)$ from this construction.

3. Let us consider the category $\catm$ whose objects are all
quadruples
$\fX=(X_*,\calps,X, \calm)$, where $X$ is a complex space, $X_*$ is
the simplicial space associated to a locally finite covering of
$X$ by Stein compact sets, $\calm$ is a complex of $\calox$-modules
with coherent cohomology bounded above, and $\calps$ is a
$\caloxs$-projective approximation, see \cite[2.19]{BFl} of the 
complex of $\caloxs$-modules
$\calms$ associated to $\calm$. With the 
morphisms in this category defined analogous to the case of
resolvents, the reader may easily
verify that this is indeed a category of models for the category $\catc$
of pairs $(X,\calm)$, where $X$ is a complex spaces and $\calm$ is a
complex of $\calox$-modules with coherent cohomology that is 
bounded above.

Using this category of models, it is then possible to define derived
symmetric and alternating powers, or even more generally arbitrary
derived Schur functors, of complexes $\calm$ as above on a complex
space $X$, and to establish the functoriality of this construction; see
also \cite{DPu} and \cite[I.4]{Ill}. Let us sketch the construction  in case of
symmetric powers.  Given a quadruple
$\fX$ as above, set
$\bbbs^p_\fX:=\calcb(\bbbs^p(\calp_*))$, where $\bbbs^p(\calps)$
denotes the
$p$-th symmetric power of the DG $\calox$-module $\calps$. 
This construction is clearly functorial in
$\fX$, and any morphism $\fX'\to \fX$ induces a quasiisomorphism
$f^{*}(\bbbs^p_\fX)\to \bbbs^p_{\fX'}$ whenever the underlying map, say
$(f,\psi)$, of pairs $(X',\calm')\to (X,\calm)$ is an isomorphism.
Applying \ref{m.5}, it follows that there are well defined derived symmetric
powers $\bbbs^p(\calm)\in D_c^-(X)$ for any $\calm\in D_c^-(X)$,
functorial in $\calm$, and even functorial with respect to
morphisms of complex spaces. 
\end{rem}

\end{document}